\documentclass[preprint]{elsarticle}
\usepackage{times,epsf,amsmath,amssymb,cite}
\usepackage{caption2}
\usepackage{hyperref}
\usepackage{amsthm}
\usepackage{graphicx} 
\usepackage{mathabx}
\usepackage{makecell}
\usepackage{wrapfig}
\usepackage[utf8]{inputenc}

\newcommand{\RNum}[1]{\uppercase\expandafter{\romannumeral #1\relax}}

\newtheorem{Theorem}{Theorem}
\newtheorem*{theorem*}{Theorem}
\newtheorem{Proposition}{Proposition}
\newtheorem{Corollary}{Corollary}
\newtheorem{Definition}{Definition}

\usepackage{bbm}
\usepackage{bm}
\newcommand\norm[1]{\left\lVert#1\right\rVert}

\newcommand\bb[1]{\mathbf{#1}}

\newcommand{\R}{\mathbb{R}}
\newcommand{\C}{\mathcal{C}}
\newcommand{\M}{\mathcal{M}}

\newcommand{\Ls}{\mathrm{L}}

\newcommand{\D}{\mathcal{D}}
\newcommand{\DN}{\mathrm{D}^N}
\newcommand{\Dd}{\mathrm{D}}
\newcommand{\K}{\mathrm{K}}

\newcommand{\w}{\mathrm{weak}^{\star}}
\newcommand{\A}{\mathrm{A}}

\journal{journal name}
\begin{document}
\begin{frontmatter}

\title{Sampling in BV-Type Spaces}

 \author[1]{Vincent Guillemet\corref{cor1}}
 \ead{vincent.guillemet@epfl.ch}
 
 \author[1]{Michael Unser}
 \ead{michael.unser@epfl.ch}

 \cortext[cor1]{Corresponding author}

 \affiliation[1]{organization={Biomedical Imaging Group,
École polytechnique fédérale de Lausanne (EPFL)},
 addressline={Station 17},
 city={Lausanne},
 postcode={CH-1015},
 state={Vaud},
 country={Suisse}}

\begin{abstract}
The sampling of functions of bounded variation (BV) is a long-standing problem in optimization. The ability to sample such functions has relevance in the field of variational inverse problems, where the standard theory fails to guarantee the mere existence of solutions when the loss functional involves samples of BV functions. 

In this paper, we prove the continuity of sampling functionals and show that the differential operator $\Dd$ admits a unique local inverse. This canonical inversion enables us to formulate an existence theorem for a class of regularized optimization problems that incorporate samples of BV functions. Finally, we characterize the solution set in terms of its extreme points.
\end{abstract}

\begin{keyword}
variational inverse problems, extreme points, approximate limits, bounded variation, splines.
\end{keyword}

\end{frontmatter}

\section{Introduction} 
Functions of bounded variations (BV) arise in the analysis of discontinuities in PDE \citep{dafermos2005hyberbolic}, in stochastic processes \citep{karatzas2014brownian,unser2014introduction}, and in image processing \citep{chambolle2016introduction,chan2005image} . They lie at the intersection of geometric measure theory and analysis \citep{evans2018measure,federer2014geometric} and extend to variational methods in free-discontinuity problems \citep{ambrosio2000functions}. They also find applications in  variational inverse problems (IPs), where one seeks to minimize data-driven functionals to recover an unknown function \citep{unser2017splines}. There, the data-acquisition modality is described by the so-called measurement functional $\bm{\nu}$.

The classic analysis of such problems first establishes the existence of minimizers 
and then characterizes the solution set. For IPs posed over BV functions and, more generally, non-reflexive Banach spaces, 
the standard theory \citep{bredies2020sparsity, flinth2019exact, unser2017splines} 
relies on the crucial assumption that $\bm{\nu}$ is $\w$-continuous. 
In many applications, however, the functional $\bm{\nu}$ 
necessarily involves sampling functionals $\delta_x:f\mapsto f(x)$. 
For general BV functions, such pointwise evaluations are not well-defined, 
except via the left and right traces $(f^-, f^+)$. 
The limitation there is that the associated sampling functionals 
\begin{equation}
    \delta_t^\pm : f \mapsto f^\pm(t)=\begin{cases}
                \underset{\epsilon\downarrow0}{\mathrm{lim }}f(t+\epsilon)\\
                \underset{\epsilon\downarrow0}{\mathrm{lim }}f(t-\epsilon)
    \end{cases}
\end{equation}
are not $\w$-continuous, a state of affairs that falls outside the limit of the standard theory 
and motivates the present analysis.

\subsection{Contributions}

Let $\Dd$ be the derivative operator, with its unique causal Green's function 
$u(\cdot)=(\cdot)_+^0$ and null space $\mathcal{N}_{\Dd}=\{1\}$ of dimension $1$. 
Since $\Dd$ has a non trivial null space, the introduction of a bi-orthogonal system $(1,\phi)$  is necessary to define a right-inverse $\Dd^{-1}_{\phi}$ \citep{unser2017splines}. 
We work in the Banach space of generalized bounded variation functions 
\begin{equation}
    \mathrm{GBV}(\R) = \bigl\{ f:\R\to\R \,:\, \norm{\Dd\{f\}}_{\M}<\infty\bigr\},
\end{equation}
equipped with the norm 
\begin{equation}
    \norm{f}_{\mathrm{GBV}}:=\norm{\Dd\{f\}}_{\M}+\vert(\mathrm{I}-\Dd^{-1}_{\phi}\circ\Dd)\{f\}\}\vert.
\end{equation}

Our main contributions are as follows.
\begin{itemize}
    \item[(i)] \textbf{Continuity of Sampling Functionals.}  
    We prove that the functional
    \begin{equation}
    \label{eq:1.2}
        \langle\delta_{t}^{\pm},\cdot\rangle:\begin{cases}
            \mathrm{GBV}(\R)\to\R\\
            f\mapsto f^{\pm}(t)
        \end{cases}       
    \end{equation}
    is continuous 
    but not $\w$-continuous.

    \item[(ii)] \textbf{Fundamental Bi-Orthogonal Systems.}  Let $\K$ be an interval. We introduce a new class of bi-orthogonal systems, that are called $\K$-\emph{fundamental systems}, 
    and show that they are unique. 
    Using the results of (i), we prove that fundamental systems yield a canonical, local inverse of $\Dd$ over $\mathrm{GBV}(\K).$

    \item[(iii)] \textbf{Existence and Structure of Variational Solutions.}  
    We apply our framework to regularized, variational inverse problems with measurement operators 
    built from functionals of the form \eqref{eq:1.2}. In Theorem \ref{th:weakRT}, we provide a nontrivial condition that guarantees that the solution set $\mathcal{V}$ is nonempty and admits a $\Dd$-spline solution. 
    Moreover, while $\mathcal{V}$ may not be $\w$-closed, we reveal that all its extreme points are $\Dd$-splines. 

    \item[(iv)] \textbf{Extension to Higher-Order BV Spaces.}  
    We extend the results of (i), (ii), and (iii) to the Banach space of functions whose $N^{\mathrm{th}}$ weak derivative is a bounded measure, and to generalized sampling functionals of the form $\Dd^{N-1}\{\delta_t^{\pm}\}.$
\end{itemize}

\subsection{Related Work}
\begin{itemize}

    \item \textbf{Inverse Problems over the Space of BV Functions.}  
    The optimization of BV functions is well established in variational IP \citep{chambolle2021approximating}. In this setting, the functional is expressed as the weighted sum of a data-fidelity term and the total-variation (TV) norm. The analysis of such problems splits into two categories. In the first one, popularized by TV denoising \citep{getreuer2012rudin}, the data are infinite-dimensional. The existence of solutions for this class of problems, as well as for the TV-flow problem, was established in \citep{andreu2001minimizing,andreu2002some,bellettini2002total}. In this case, the solution does not admit a sparse representation.
    
    \item \textbf{Sparsity Without Sampling Functionals.}  
    If the data are finite-dimensional and belong to $\R^M$, standard results \citep{boyer2019representer,bredies2020sparsity,flinth2019exact,unser2017splines} show that the extreme points of the solution set are $\Dd$-splines with at most $M$ knots. This sparsity result holds as long as the measurement functional is $\w$-continuous, for example when it is defined in terms of indicator functionals $\mathbbm{1}_{A}$. Nevertheless, when the measurement operator involves non-$\w$-continuous left and right traces, the theory breaks down and fails to even ensure the mere existence of a solution. One possible remedy is to use a regularization term involving higher-order derivatives, which promotes candidates with greater smoothness. Relevant examples include the generalized total variation (GTV) \citep{bredies2020sparsity,bredies2014regularization,bredies2010total} and the Hessian total variation (HTV) \citep{ambrosio2024linear,lefkimmiatis2011hessian,pourya2024box}. The trade-off, however, is that the corresponding function spaces are strictly smaller than BV, and therefore cannot represent discontinuities.
    
    \item \textbf{Sampling Functionals Without   Sparsity.}  
    A fundamental family of problems involving left and right traces is referred to as \emph{free discontinuity problems} by E.~De Giorgi \citep{de1991free}. In such problems, the functional is defined as the combination of a volume energy term and a surface energy term. In dimension~1, the surface energy term simply involves samples of BV functions. The authors of \citep{ambrosio2000functions} provide an existence result for a broad class of free-discontinuity problems. In particular, their result applies to the Mumford--Shah functional \citep{mumford1989optimal}, which has been a cornerstone of image processing \citep{morel2012variational}. Despite their generality, these results are not applicable within the variational IP framework. 

\end{itemize}

\section{Mathematical Background}

\subsection{Continuity in Dual Banach Spaces}
Let $\mathcal{B}, \mathcal{B}^{\star}$, and $\mathcal{B}^{\star\star}$ be three Banach spaces such that 
$\big(\mathcal{B}^{\star},\norm{\cdot}_{\mathcal{B}^{\star}}\big)$ is the dual space of $\big(\mathcal{B},\norm{\cdot}_{\mathcal{B}}\big)$ and $\big(\mathcal{B}^{\star\star},\norm{\cdot}_{\mathcal{B}^{\star\star}}\big)$ is the dual space of $\big(\mathcal{B}^{\star},\norm{\cdot}_{\mathcal{B}^{\star}}\big)$.  
A sequence $(f_n)_{n=1}^{\infty}$ of elements $f_n\in\mathcal{B}^{\star}$ is said to be $\w$-convergent to a limit $f\in\mathcal{B}^{\star}$ if 
\begin{equation}
    \forall g\in\mathcal{B}:\quad\lim_{n\to\infty}\,\langle f_n,g\rangle=\langle f,g\rangle.
\end{equation}
A linear functional 
\begin{equation}
    \eta:\mathcal{B}'\to\R
\end{equation}
\begin{itemize}
    \item is continuous if, for every sequence $(f_n)_{n=1}^{\infty}$ converging in norm to some $f$, one has that
    \begin{equation}
        \lim_{n\to\infty}\eta(f_n)=\eta(f);
    \end{equation}
    moreover, the continuity of $\eta$ is equivalent to the existence of a representative $h\in\mathcal{B}^{\star\star}$ such that 
    \begin{equation}
        \forall f\in\mathcal{B}^{\star}:\quad \eta(f)=\langle h,f\rangle;
    \end{equation}
    \item is $\w$-continuous if, for every sequence $(f_n)_{n=1}^{\infty}$ $\w$-convergent to some $f$, one has that
    \begin{equation}
        \lim_{n\to\infty}\eta(f_n)=\eta(f);
    \end{equation}
    moreover, the $\w$-continuity of $\eta$ is equivalent to the existence of a representative $h\in\mathcal{B}$ such that 
    \begin{equation}
        \forall f\in\mathcal{B}^{\star}:\quad \eta(f)=\langle f,h\rangle.
    \end{equation}
\end{itemize}
Since $\mathcal{B}\subset\mathcal{B}^{\star\star}$, if $\eta$ is $\w$-continuous, then it is also continuous. The converse is false for infinite-dimensional (non-reflexive) vector spaces. Therefore, contrarily to what the terminology suggests, continuity is weaker than $\w$-continuity.

\subsection{Native Space-BV Theory}
\label{sec:2.2}
To review the key results of \citep{unser2019native,unser2017splines}, it is necessary to introduce the notion of admissible systems $(1, \phi)$ and their associated kernel $g_{\phi}$. In the sequel, we denote the space $\mathrm{GBV}(\R)$ by $\M_{\Dd}(\R)$.
\begin{Definition}
\label{def:admissible}
The system $(1, \phi)$ is a universal $\Dd$-admissible system if the function $\phi\in\mathcal{D}(\R)$ with $\mathrm{supp}(\phi)\subset[\phi^-,\phi^+]$ satisfies 
    \begin{equation}
    \label{eq:3.13}
            \langle 
            1, \phi\rangle=1.
    \end{equation}
\end{Definition}
The kernel $g_{\phi},$ and the associated integral transform $\Dd_{\phi}^{-1}$ are defined as (see \citep{unser2017splines}) 
\begin{equation}
\label{eq:2.8}
  g_{\phi}:(t,x)\mapsto g_{\phi}(t,x)=u(t-x)-\langle\phi(\cdot), u(\cdot-x)\rangle,
\end{equation}
\begin{equation}
    \Dd_{\phi}^{-1}:f\mapsto\int_{\R}g_{\phi}(t,x)f(x)\mathrm{d}x.
\end{equation}
We recall Theorem \ref{th:1}, which establishes several key properties of $\M_\Dd(\R)$. The projection operator $\text{proj}_{\mathcal{N}_{\Dd}}$ is defined as 
\begin{equation}
    \text{proj}_{\mathcal{N}_{\Dd}}:
    \begin{cases}
        \D'(\R)\to\mathcal{N}_\Dd\\
        f\mapsto\langle \phi,f\rangle.
    \end{cases}
\end{equation}
\begin{Theorem}[\citep{unser2017splines}]
\label{th:1}
If $(1, \phi)$ is a universal $\Dd$-admissible system, then the space $\M_{\D}(\R)$ admits the direct-sum decomposition
\begin{equation}
\label{eq:3.23}
    \M_{\Dd}(\R)=\Dd_{\phi}^{-1}\{\M(\R)\}\oplus\mathcal{N}_{\Dd},
\end{equation}
with the norm 
\begin{equation}
    \norm{\cdot}_{\M_{\Dd}}:=\norm{\Dd\{\cdot\}}_{\M}+\norm{\mathrm{proj}_{\mathcal{N}_{\Dd}}\{\cdot\}}_{1},
\end{equation}
for which it is a Banach space. 
\end{Theorem}

A predual is given by 
$\C_{\Dd}(\R)=\Dd^{\star}\{\C_0(\R)\}\oplus\mathcal{N}_N^{\star}$ \citep{unser2017splines}. This characterization enables us to extend the notion of an admissible system to $\phi\in\C_{\Dd}(\R)$ \citep{unser2019native}. Indeed, we observe \emph{a posteriori} that the construction of the space $\M_{N}(\R)$, as well as Theorem \ref{th:1}, hold for any $\Dd$-admissible system $(1, \phi),$ with all norms being equivalent.

An important class of functions in $\M_{\Dd}(\R)$ is the class of $\Dd$-splines.
\begin{Definition}
\label{def:spline}
    A distribution $f\in\D'(\R)$ is a $\Dd$-spline if 
    \begin{equation}
        \Dd\{f\}=\sum_{k=1}^Ka_k\delta_{x_k}
    \end{equation}
    for some $K\in\mathbb{N}$, coefficients $(a_k)_{k=1}^K\in\R^K$, and knots $(x_k)_{k=1}^K\in\R^K$. If we define the innovation $\mu
    =\sum_{k=1}^Ka_k\delta_{x_k}\in\D'(\R)$, then we can generate such a $\Dd$-spline as
    \begin{equation}
        f=c_0+(u\ast \mu),\quad\text{with}\quad c_0\in\mathcal{N}_{\Dd},
    \end{equation}
which is such that 
\begin{equation}
    \norm{f}_{\M_{\Dd}}=\vert c_0\vert+\sum_{k=1}^K\vert a_k\vert.
\end{equation}
\end{Definition}

\section{Sampling}
\label{sec:3}
We introduce the left- and right-limit functionals $\delta_t^-$ and $\delta_t^+$, respectively, such that, for all $f\in\M_{\Dd}(\R)$,
\begin{equation}
    \langle\delta_{t}^-,f\rangle=\underset{\epsilon\downarrow0}{\mathrm{lim }}f(t-\epsilon),\quad\text{and}\quad\langle\delta_{t}^{+},f\rangle=\underset{\epsilon\downarrow0}{\mathrm{lim }}f(t+\epsilon),
\end{equation}
and further define their convex combination, for all $\omega\in[0,1]$,
\begin{equation}
  \delta_{t}^{\omega}=\omega\delta_{t}^{+}+(1-\omega)\delta_{t}^{-}.
\end{equation}
Observe that the sampling functional
\begin{equation}
\label{eq:2.2.22}
\langle\delta_{t}^{\omega},\cdot\rangle:\begin{cases}
    \M_{\Dd}(\R)\to\R\\
    f\mapsto\langle\delta_{t}^{\omega},f\rangle
\end{cases}       
\end{equation}
is not $\w$-continuous. Indeed, for $\omega=1$ and $t=0$, the sequence $(f_\ell)_{\ell=1}^{\infty}=(u\big(\cdot-\frac{1}{\ell}\big))_{\ell=1}^{\infty}$ is $\w$-convergent to $u$, but 
\begin{equation}
    \lim_{\ell\to\infty}\langle\delta_0^{+},f_\ell\rangle= \lim_{\ell\to\infty}\left\langle\delta_0^{+},\left(\cdot-\frac{1}{\ell}\right)_+^0\right\rangle=0\neq\langle\delta_0^{+},u\rangle=1.
\end{equation}
Although \eqref{eq:2.2.22} is not $\w$-continuous, we state in Proposition \ref{diractweak} that it is continuous.
\begin{Proposition}
\label{diractweak}
The mapping in \eqref{eq:2.2.22} is well-defined, linear, and continuous. 
\end{Proposition}

\begin{proof}[\textbf{Proof of Proposition \ref{diractweak}.}]
It is sufficient to prove the stated properties for the left and right limits. We only provide a proof for the right limit. Equation \eqref{eq:2.2.22} is clearly well-defined, linear, and continuous over $\mathcal{N}_{\Dd}$. By invoking Theorem \ref{th:1}, we fix $f\in\Dd_{\phi}^{-1} \{\M(\R)\}$. Next, we recall from \citep{guillemet2025adaptive} that
\begin{equation}
\label{eq:supp}
    \forall t\in\R,\forall x\notin[\mathrm{min}(t,\phi^-),\mathrm{max}(t,\phi^+)]:\quad g_{\phi}(t,x)=0.
\end{equation}
Let $\A_{\epsilon}=[\text{min}(t,\phi^-),\text{max}(t+\epsilon,\phi^+)]$. Then, the combination of \eqref{eq:2.8} and \eqref{eq:supp} yields, for all $\tau\in[t,t+\epsilon]$,
\begin{align}
    f(\tau)&=\int_{\A_{\epsilon}}u(\tau-x)\,\mathrm{d}\mu(x)-\int_{\A_{\epsilon}}\langle \phi(\cdot),u(\cdot-x)\rangle\,\mathrm{d}\mu(x)\nonumber\\
    &=\mu([\mathrm{min}(t,\phi^-),\tau])-\int_{\A_{\epsilon}}\langle \phi(\cdot),u(\cdot-x)\rangle\,\mathrm{d}\mu(x)\nonumber\\
    \Rightarrow\,\lim_{\epsilon\downarrow0}f(t+\epsilon)&=\mu([\text{min}(t,\phi^-),t])-\int_{\A_0}\langle \phi(\cdot),u(\cdot-x)\rangle\,\mathrm{d}\mu(x).
\end{align}
Consequently, the right limit is well-defined, linear, and its continuity follows from the upper bound
\begin{equation}
    \vert\langle\delta_{t}^{+},\Dd^{-1}_{\phi}\{\mu\}\rangle\vert\leq\norm{\mu}_{\M}\left(\underset{x\in A_0}{\text{max}}\vert\langle \phi(\cdot),u(\cdot-x)\rangle\vert+1\right).
\end{equation}
\end{proof}
\noindent Therefore, for $\M_{\Dd}(\R)^{\star}$, which is the dual space of $\M_{\Dd}(\R)$, one has that
\begin{equation}
\delta_{t}^{\omega}\notin\C_{\Dd}(\R),\quad\text{but}\quad\delta_{t}^{\omega}\in\M_{\Dd}(\R)^{\star}.
\end{equation}


\section{Fundamental Systems}
As demonstrated in Section \ref{sec:2.2}, the inversion of $\Dd$ on $\M(\R)$ involves the use of a non-shift-invariant kernel $g_{\phi}$, itself tailored to a $\Dd$-admissible system $(1,\phi)$. Leveraging our understanding of sampling functionals, we establish that, under suitable constraints, the choice of $(1,\phi)$ and $g_{\phi}$ is unique.

As it may be inconvenient to work directly in the space $\M_{\Dd}(\R)$, the first constraint is about localization. In fact, in some applications, one is interested in the value of $f\in\M_{\Dd}(\R)$ only on a known interval $\mathrm{K}=[K^-,K^+].$ We now specialize the notion of $\Dd$-admissible system to take this observation into account

\subsection{Localization}
\begin{Definition}
\label{def:localizedsystem}
    A $\Dd$-admissible system $(1,\phi)$ is called $\K$-localized if $[\phi^-,\phi^+]\subset[K^-,K^+]$.
\end{Definition}
In Proposition \ref{prop:localizedjustification}, we show that $\K$-localized systems yield a favorable weight distribution of the kernel $g_{\phi}$.

\begin{Proposition}
\label{prop:localizedjustification}
    If $(1,\phi)$ is a $\Dd$-admissible system, then 
\begin{equation}
\label{eq:loc.1}
    \forall x\notin[\mathrm{min}(\phi^-,K^-),\mathrm{max}(\phi^+,K^+)],\ \forall t\in[K^-,K^+]:\quad g_{\phi}(t,x)=0.
\end{equation}
    If $(1,\phi)$ is a $\K$-localized system, then 
\begin{equation}
\label{eq:loc.2}
    \forall x\notin[K^-,K^+],\ \forall t\in[K^-,K^+]:\quad g_{\phi}(t,x)=0
\end{equation}
and $\forall f\in\M_{\Dd}(\R),\ \forall t\in[K^-,K^+],$
\begin{equation}
\label{eq:loc.3}
    f(t)=c_0+\int_{\R}g_{\phi}(t,x)\mathrm{d}\mu(x)=c_0+\int_{\K}g_{\phi}(t,x)\mathrm{d}\mu(x).
\end{equation}
\end{Proposition}

\begin{figure}
    \centering
    \includegraphics[width=1\linewidth]{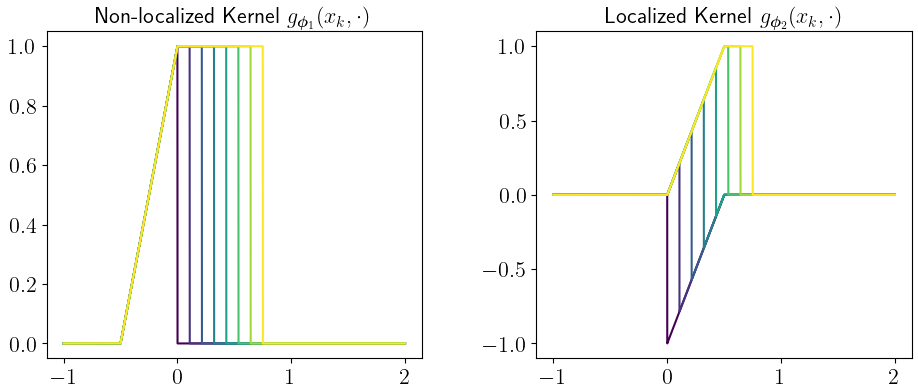}
    \caption{The kernel $g_{\phi_1}$ is constructed from the non-$[0,1]$-localized system $(1,2\mathbbm{1}_{[-0.5,0]})$. The kernel $g_{\phi_2}$ is constructed from the $[0,1]$-localized system $(1,2\mathbbm{1}_{[0,0.5]})$. The color of the graph reflects the second dimension of the kernels. The darker the graph, the smaller the value of the corresponding knot $x_k$, with the six knots $x_k$ being uniformly sampled in $[0,0.75]$.}
    \label{fig:1}
\end{figure}

\begin{proof}[\textbf{Proof of Proposition \ref{prop:localizedjustification}}]
 Recall that the kernel $g_{\phi}(t,x)$ is such that 
\begin{equation}
       \begin{cases}
           \forall x<\phi^-,\ \forall t>x\\
           \forall x>\phi^+,\ \forall t<x
        \end{cases}:\quad g_{\phi}(t,x)=0.
   \end{equation}
Therefore,
\begin{itemize}
    \item if $\phi^-\leq K^-$, then $\forall x<\phi^-$ it holds that $\forall t\geq K^-\geq\phi^->x:\ g_{\phi}(t,x)=0;$
    \item if $\phi^->K^-$, then $\forall x<K^-$ it holds that $\forall t\geq K^->x:\ g_{\phi}(t,x)=0;$
    \item if $\phi^+\geq K^+$, then $\forall x>\phi^+$ it holds that $\forall t\leq K^+\leq\phi^+<x:\ g_{\phi}(t,x)=0;$
    \item if $\phi^+<K^+$, then $\forall x>K^+$ it holds that $\forall t\leq K^+<x:\ g_{\phi}(t,x)=0.$
\end{itemize}
Consequently,
\begin{equation}
\label{app:A.8.120}
    \forall x\notin[\text{min}(\phi^-,K^-),\text{max}(\phi^+,K^+)],\ \forall t\in[K^-,K^+]:\quad g_{\phi}(t,x)=0.
\end{equation}
This shows \eqref{eq:loc.1}, while \eqref{eq:loc.2} and \eqref{eq:loc.3} follow directly.
\end{proof}

A stronger constraint would be to require that \eqref{eq:loc.3} holds $\forall t\in\R$. To clarify its implications, we recall that the definition of the Banach subspace $\M_{\Dd}(\K)$ \citep{guillemet2025convergence} is
\begin{equation}
\label{eq:locnative}
    \M_{\Dd}(\K)=\{f\in\M_{\Dd}(\R):\quad\Dd\{f\}\in\M(\K)\}.
\end{equation}
In contrast with functions in $\M_{\Dd}(\R)$, those in $\M_{\Dd}(\K)$ admit an integral representation with a shift-invariant kernel. 
\begin{Proposition}[\citep{guillemet2025convergence}]
\label{prop:6}
    A distribution $f\in\mathcal{D}'(\R)$ belongs to the space $\M_{\Dd}(\K)$ if and only if it is of the form, $\forall t\in\R,$
\begin{align}
    f(t)& =\tilde{c}_{0}+\int_{\K}g_{\phi}(t,x)\mathrm{d}\mu(x),\label{eq:4.36}\\
    &=c_0+\int_{\K}u(t-x)\mathrm{d}\mu(x),\label{eq:fund.33}
\end{align}
with $\mu\in\M(\K)$ and $c_0,\tilde c_0\in\R.$
\end{Proposition} 
This representation allows us to provide a new characterization of the action of the left- and right-limit functionals on $\M_{\Dd}(\K)$. In particular, Proposition \ref{prop:weakrep} shows that the left and right limits at $t$ differ if and only if the supporting measure $\mu$ of $f$ has an atom at $t$.
\begin{Proposition}
    \label{prop:weakrep}
For all $f\in\M_{\Dd}(\K)$ parameterized by $(\mu,c_0)\in\M(\K)\times\R$ as in \eqref{eq:fund.33}, and $\forall t\in\K$,  one has that 
\begin{equation}
\label{eq:3.1.37}
\langle \delta_{t}^{-},f\rangle=c_0+\mu([K^-,t[)
\end{equation}
and
\begin{equation}
\label{eq:3.1.38}
\langle\delta_{t}^{+},f\rangle=c_0+\mu([K^-,t]).
\end{equation}
\end{Proposition}

\begin{proof}[\textbf{Proof of Proposition \ref{prop:weakrep}}]
We prove only \eqref{eq:3.1.37}, by stating that
\begin{align}
    \langle\delta_{t}^{-},u\ast \mu\rangle=\lim_{\epsilon\downarrow0}\int_{\K}u(t- \epsilon-x)\,\mathrm{d}\mu(x)&=\lim_{\epsilon\downarrow0}\mu([K^-,t-\epsilon[).
\end{align}
\end{proof}
Unlike \eqref{eq:4.36}, and owing to the non–shift-invariant structure of $g_{\phi}$, the sum in the representation \eqref{eq:fund.33} is not direct. For the Banach space $\M_{\Dd}(\R)$, because the domain $\R$ is unbounded, it is impossible to let $g_{\phi}$ specified by \eqref{eq:2.8} be shift-invariant. In addition to the localization constraint described in this section, we impose in Section \ref{sec:3.2} a new, \emph{local} shift-invariance constraint. Its enforcement is made possible by specific properties of $\M_{\Dd}(\K)$. The key is provided in Proposition \ref{diracweakstar2}, where we build on Proposition \ref{prop:weakrep} to show that certain left- and right-limit functionals are $\w$-continuous over $\M_{\Dd}(\K)$. First, we establish in Proposition \ref{prop:predual} that $\M_{\Dd}(\K)$ is a dual Banach space whose $\w$-topology coincides with the restricted $\w$-topology from $\M_{\Dd}(\R)$.

\begin{Proposition}
\label{prop:predual}
Define the annihilator
\begin{equation}
^{\perp}\M_{\Dd}(\K) := \{v\in\C_{\Dd}(\R) : \forall f\in\M_{\Dd}(\K),\quad\langle f,v\rangle=0 \}.
\end{equation}
Then, there exists a canonical isometric isomorphism $I$ such that  
\begin{equation}
    \M_{\Dd}(\K)\overset{I}{\cong}\big(\C_{\Dd}(\R)/^{\perp}\M_{\Dd}(\K)\big)^{\star}.
\end{equation}
In particular, $\C_{\Dd}(\R)/^{\perp}\M_{\Dd}(\K)$ can be identified as the predual of $\M_{\Dd}(\K).$
\end{Proposition}

\begin{proof}[\textbf{Proof of Proposition \ref{prop:predual}}]
We define the quotient map
\begin{equation}
    q:\C_{\Dd}(\R)\to\C_{\Dd}(\R)/^{\perp}\M_{\Dd}(\K),
\end{equation}
together with the identification map 
\begin{equation}
I:\M_{\Dd}(\K)\mapsto\big(\C_{\Dd}(\R)/^{\perp}
\M_{\Dd}(\K)\big)^{\star},
\end{equation}
such that 
\begin{equation}
\forall[v]\in\C_{\Dd}(\R)/^{\perp}\M_{\Dd}(\K):\quad I(f)([v]) = \langle f,v\rangle.  
\end{equation}
It is straightforward to verify that $I$ is well-defined, linear, and injective. 

We now prove surjectivity. Let $\Phi\in\big(\C_{\Dd}(\R)/^{\perp}\M_{\Dd}(\K)\big)^{\star}$ and define $f_{\Phi}=\Phi\circ q\in\M_{\Dd}(\R).$ Then, $f_{\Phi}$ vanishes on $^{\perp}\M_{\Dd}(\K)$, with
\begin{equation}
    f_{\Phi}\in\Big(^{\perp}\M_{\Dd}(\K)\Big)^{\perp}:=\big\{f\in\M_{\Dd}(\R):\quad\forall v\in^{\perp}\M_{\Dd}(\K),\quad\langle f,v\rangle=0\big\}.
\end{equation}

Since $\M_{\Dd}(\K)$ is $\w$-closed, it follows from \citep[Theorem 4.7]{rudinfunction} that $\Big(^{\perp}\M_{\Dd}(\K)\Big)^{\perp}=\M_{\Dd}(\K).$ Hence $f_{\Phi}\in\M_{\Dd}(\K)$ and $\forall[v]\in\C_{\Dd}(\R)/^{\perp}
\M_{\Dd}(\K),$ 
\begin{equation}
    I(f_{\Phi})([v])=\langle f_{\Phi},v\rangle=\Phi(q(v))=\Phi([v]).
\end{equation}
Therefore, $I(f_{\Phi})=\Phi$, and $I$ is surjective. It remains to show that $I$ is an isometry. Denote by $\norm{\cdot}^{'}$ the norm of $\big(\C_{\Dd}(\R)/^{\perp}
\M_{\Dd}(\K)\big)^{\star}$. We  observe that, $\forall f\in\M_{\Dd}(\K),$
\begin{align}
    \norm{I(f)}^{'}&=\underset{[v]:\underset{w\in[v]}{\mathrm{inf}}\norm{w}_{\C_{\Dd}(\R)}\leq1}{\mathrm{sup}}I(f)([v])=\underset{v:\underset{w\in[0]}{\mathrm{inf}}\norm{w+v}_{\C_{\Dd}(\R)}\leq1}{\mathrm{sup}}\langle f,v\rangle\nonumber\\
    &=\underset{v:\norm{v}_{\C_{\Dd}(\R)}\leq1}{\mathrm{sup}}\langle f,v\rangle=\norm{f}_{\M_{\Dd}(\R)}.
\end{align}
\end{proof}
Proposition \ref{prop:predual} establishes that $\M_{\Dd}(\K)$ also admits a $\w$-topology, such that $(f_n)_{n=1}^{\infty}$ is $\w$-convergent to $f$ in $\M_{\Dd}(\K)$ if and only if
\begin{equation}
    \forall [v]\in\C_{\Dd}(\R)/^{\perp}\M_{\Dd}(\R):\quad\lim_{n\to\infty}\langle f_n,[v]\rangle=\langle f,[v]\rangle.
\end{equation}
Note that, although the quotient appears to complicate the predual structure, one still has that $(f_n)_{n=1}^{\infty}$ is $\w$-convergent to $f$ if and only if
\begin{equation}
    \forall v\in\C_{\Dd}(\R):\quad\lim_{n\to\infty}\langle f_n,v\rangle=\langle f,v\rangle,
\end{equation}
because a change in the representative in the quotient class does not alter the value of the brackets.

\begin{Proposition}
\label{diracweakstar2}
The mappings
\begin{equation}
\label{eq:3.1.41}
\langle\delta_{K^{\pm}}^{\pm},\cdot\rangle:\begin{cases}
    \M_{\Dd}(\K)\to\R\\
    f\mapsto\langle\delta_{K^{\pm}}^{\pm},f\rangle=f^{\pm}(K^{\pm})
\end{cases}       
\end{equation}
are $\w$-continuous.
\end{Proposition}

\begin{proof}[\textbf{Proof of Proposition \ref{diracweakstar2}}]
We only prove the $\w$ continuity of the left limit. Assume that $f$ is as in \eqref{eq:fund.33}. Then, by Proposition \ref{prop:weakrep} and for $\K_T=[K^--T,K^-],$
\begin{align}
    \langle\delta_{t}^{-},f\rangle=\left\langle\frac{\mathrm{1}_{\K_T}(\cdot)}{T},f
    \right\rangle=\left\langle\frac{\mathrm{1}_{K_T}(\cdot)}{T},c_0
    \right\rangle=\int_{\R}\frac{\mathrm{1}_{K_T}(x)}{T}c_0\,\mathrm{d}x=c_0.
\end{align}
Consequently, the mapping $\langle\delta_{K^{-}}^{-},\cdot\rangle$ can be represented by the functional $\frac{\mathrm{1}_{\K_T}(\cdot)}{T}$, which belongs to $\C_{\Dd}(\R)$ by \citep[Corollary 1]{guillemet2025adaptive}. The $\w$ continuity follows directly.
\end{proof}

The subtlety of Proposition \ref{diracweakstar2} lies in the fact that, although the functionals in \eqref{eq:3.1.41} are merely continuous over $\M_{\Dd}(\R)$, they are $\w$-continuous over $\M_{\Dd}(\K)$, owing to the restriction that functions in $\M_{\Dd}(\K)$ may be discontinuous only on $\K$.

\subsection{Locally Shift-Invariant Kernels}
\label{sec:3.2}

We introduce in Definition \ref{def:idealsystem} the local shift-invariance constraint for localized systems.

\begin{Definition}
\label{def:idealsystem}
    A $\K$-localized system $(1,\phi)$ is called $\K$-fundamental if
    \begin{itemize}
        \item[1.] the functional $\phi:\M_{\Dd}(\K)\to\R$ is $\w$-continuous;
        \item[2.] the kernel $g_{\phi}$ reduces to the Green’s function, in the sense that
    \begin{equation}
    \forall x\in[K^-,\infty[:\quad g_{\phi}(\cdot,x)=u(\cdot-x).
    \end{equation}
    \end{itemize}
\end{Definition}

In contrast with $(1,\phi)$, if a system is $\K$-fundamental, then it is denoted by $(1,\iota)$. A distinctive feature is that a $\K$-fundamental system is only required to be $\w$-continuous on $\M_{\Dd}(\K)$ (Definition \ref{def:idealsystem}, Item 1), and not on $\M_{\Dd}(\R)$. This choice of kernel leads to a significant simplification of the theory, as it yields the representation \eqref{eq:fund.33} while it preserves the direct-sum structure. Indeed,
\begin{equation}
\label{eq:fund.36}
    \M_{\Dd}(\K)=\Dd_{\iota}^{-1}\{\M(\K)\}\oplus\mathcal{N}_{\Dd}=\Dd^{-1}\{\M(\K)\}\oplus\mathcal{N}_{\Dd},
\end{equation}

\begin{figure}[h!]
    \centering
    \includegraphics[width=1\linewidth]{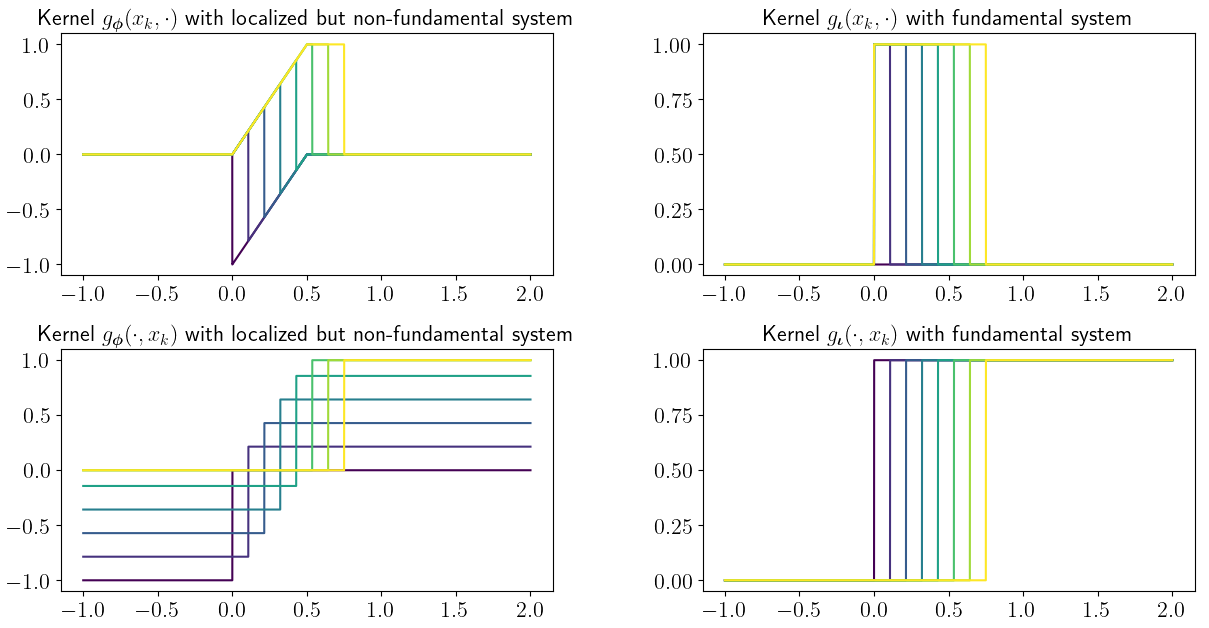}
    \caption{The kernel $g_{\phi}$ is constructed from the $[0,1]$-localized but non-fundamental system $(1,2\mathbbm{1}_{[0,0.5]})$. The first column presents graphs of the kernel $g_{\phi}$, while the second column presents graphs of $g_{\iota}$. In the first row, the color of the graph reflects the second dimension of the kernels, whereas in the second row it reflects the first dimension. The darker the graph, the smaller the value of the corresponding knot $x_k$, with the six knots $x_k$ being uniformly sampled in $[0,0.75]$.
}
    \label{fig:2}
\end{figure}

\noindent where the integral operator $\Dd^{-1}$ is described by the kernel $u$ and 
\begin{equation}
    \Dd^{-1}\{\M(\K)\}=\left\{\int_{\K}u(\cdot-x)\,\mathrm{d}m(x):\ m\in\M(\K)\right\}.\label{eq:4.4.76}
\end{equation}
In \eqref{eq:fund.36}, the first equality follows from \eqref{eq:3.23} and the definition of $\M_{\Dd}(\K)$, while the second equality follows from the fact that $(1,\iota)$ is $\K$-fundamental. We now show that the $\K$-fundamental system is unique.

\begin{Theorem}
 \label{prop:ideal2}
The system $(1,\iota)=(1,\delta_{K^-}^-)$ is the unique $\K$-fundamental system.
\end{Theorem}

\begin{proof}[\textbf{Proof of Theorem \ref{prop:ideal2}.}]\hfill\\
\textbf{$\K$-localized and Fundamental.} 
The property that $(1,\delta_{K^-}^-)$ is $\K$-localized is straightforward to verify. Item 1 of the $\K$-fundamental definition is ensured by  Proposition \ref{diracweakstar2}. Item~2 is verified as follows: $\forall x\in[K^-,\infty[,\,\forall t\in\R$,
\begin{equation}
    g_{\phi}(t,x)=u(t-x)-\langle\iota, u(\cdot-x)\rangle,
\end{equation}
with
\begin{align}
   \langle\iota, u(\cdot-x)\rangle=\lim_{\epsilon\downarrow0}u(K^--x-\epsilon)=0.
\end{align}
\textbf{Uniqueness.} Assume that $(1,\tilde{\iota})$ is another $\K$-fundamental system. Observe that, with respect to 
Proposition \ref{prop:predual}, any $\w$-continuous linear map over $\M_{\Dd}(\K)$ is in fact represented by an equivalence class. We denote this class by $[\tilde \iota]$. There exists a representative $\hat{\iota}$ of the quotient class $[\tilde{\iota}]$ such that 
\begin{equation}
    \hat{\iota}\in\C_{\Dd}(\R)=\Dd^{\star}\{\C_0(\R)\}\oplus\text{span}\{\phi\},\quad\text{with}\quad\phi\in\D(\R).
\end{equation}
Consider the representation $\hat{\iota}=\Dd^{\star}\{v\}+a\phi$. Item~2 of the fundamental assumption implies that, $\forall x\in]K^-,\infty[,$
\begin{align}
    g_{\tilde{\iota}}(\cdot,x)=u(\cdot-x)&\Leftrightarrow\langle\tilde{\iota}_{n}, u(\cdot-x)\rangle=0\nonumber\\
    &\Leftrightarrow\langle\hat{\iota}_{n}, u(\cdot-x)\rangle=0\nonumber\\
    &\Leftrightarrow v(x)=-\left\langle a\phi, u(\cdot-x)\right\rangle\nonumber\\
    &\Leftrightarrow v(x)=-\left(u^{\vee}\ast a\phi\right)(x).\label{eq:A.6.107}
\end{align}
By introducing \eqref{eq:A.6.107} into the representation of $\hat{\iota}_{n}$, we find that, $\forall\rho\in\D(\R)$ such that $\text{supp}\{\rho\}\subset]K^-,\infty[,$
\begin{equation}
    \langle\hat{\iota}_{n},\rho\rangle=\langle\Dd^{\star}\{-u^{\vee}\ast a\phi\}+a\phi,\rho\rangle=\langle-a\phi+a\phi,\rho\rangle=0.
\end{equation}
It follows that $\text{supp}\{\hat{\iota}\}\subset]-\infty,K^-]$ and, therefore, that the support of any representative 
in $[\tilde{\iota}]$ must be a subset of $]-\infty,K^-]$. In particular, $\tilde{\iota}$ must satisfy $\text{supp}\{\tilde{\iota}\}=]-\infty,K^-]\cap[K^-,\infty[=\{K^-\}$. The space of admissible functionals that are $\w$-continuous over $\M_{\Dd}(\K)$, and whose support equals $\{K^-\}$ is spanned by the functional $\delta_{K^-}^-.$ This proves the claim.
\end{proof}

Theorem \ref{prop:ideal2} reveals that, over $\M_{\Dd}([0,T])$, where $[K^-,K^+]=[0,T]$, the unique appropriate system is $(1,\delta_0^-)$ which, in turn, establishes that the unique appropriate norm is
\begin{equation}
    \norm{\Dd\{f\}}_{\M}+\vert f(0^-)\vert
\end{equation}
which induces a direct-sum topology with respect to the kernel $u(\cdot).$ We now demonstrate that the extreme-point solutions of regularized optimization problems, with a regularization term on the null space, are $\Dd$-splines.

\begin{Corollary}
\label{RTnull space}
If $\bm{\nu}=(\nu_m)_{m=1}^M\in\C_{\Dd}(\R)^{M}$ and $\bb{y}\in\R^M$, then the solution set 
\begin{equation}
    \mathcal{V}=\underset{f\in\M_{\Dd}(\K)}{\mathrm{argmin}}\,\Big(\norm{\bb{y}-\langle\bm{\nu},f\rangle}_2^2+\lambda\norm{f}_{\M_{\Dd}}\Big)
\end{equation}
is nonempty, $\w$-compact, and convex. Its extreme points are of the form 
\begin{equation}
\label{eq:3.2.58}
    \sum_{k=1}^{K}a_kg_{\phi}(\cdot,x_k)+c,\quad (K+\norm{c}_{0})\leq M,\quad\text{where}\quad\norm{c}_0=\begin{cases}
        0,&c=0\\
        1,&\text{else}.
    \end{cases}
\end{equation}
Moreover, if the underlying system is the unique $\K$-fundamental system, then the extreme points are of the form
\begin{equation}
\label{eq:3.2.59}
    \sum_{k=1}^{K}a_ku(\cdot-x_k)+c,\quad (K+\norm{c}_{0})\leq M.
\end{equation}
\end{Corollary}

The proof of \eqref{eq:3.2.58} follows the standard reasoning in \citep{guillemet2025convergence,unser2017splines} for representer theorems, while the proof of \eqref{eq:3.2.59} follows from Item~2 of the assumption that the system is fundamental. Corollary \ref{RTnull space} shows that, in contrast with the regularization-free null-space scenario, once regularization is introduced, the extreme points are represented sparsely with respect to a shift-invariant kernel if and only if the underlying system is fundamental.

\subsection{Causality}
So far, we have worked with the (unique) causal Green's function $u$, and we have also constructed localized systems, showing that there exists a unique \emph{causal} fundamental system. It is interesting to remark that a similar construction can be carried out with the (unique) anti-causal Green's function $-u(-\cdot)=-u^{\vee}(\cdot)$. In this case, the unique anti-causal fundamental system $(1,\iota)$ is given by $(1,\delta_{K^+}^+)$. It follows from its existence that every function $f\in\M_{\Dd}(\K)$ admits the direct-sum decomposition
\begin{equation}
    f=c+(u^{\vee}\ast m),\quad m\in\M(\K),\; c\in\mathcal{N}_{\Dd}.
\end{equation}
In what follows, fundamental systems are employed in their causal version, unless specified otherwise.


\section{Resolution of Continuous-Domain Optimization Problems}
\subsection{Definitions and Assumptions}
Let $\bb{t}=\{t_m\}_{m=1}^M$ and $\K=[K^-,K^+]$ be such that $t_m\in\K.$ Let $\mathrm{X}\subset\K$ be another compact set. We are interested in optimization problems as in
\begin{equation}
\label{eq:3.1.30}
\mathcal{V}(\mathrm{X})=\underset{f\in\M_{\Dd}(\mathrm{X})}{\text{argmin}}\quad\mathcal{J}(f),
\end{equation}
with the $\emph{loss functional}$
\begin{equation}
\label{eq:3.1.31}
\mathcal{J}(f)=E\left(\bb{y},\left\langle\begin{bmatrix}
\bm{\bm{\delta}_{\bb{t}}}^{\omega}\\
\bm{\nu}
\end{bmatrix},f\right\rangle\right)+\lambda\norm{f}_{\M_{\Dd}},
\end{equation}
where $\bb y\in\R^{M+\tilde M}$, $E$ is the data-fidelity functional, and where 
\begin{equation}
   \langle\bm{\bm{\delta}_{\bb{t}}}^{\omega},f\rangle=\Big(\langle\delta_{t_m}^{\omega},f\rangle\Big)_{m=1}^M,\quad\quad\langle\bm{\nu},f\rangle=\Big(\langle\nu_m,f\rangle\Big)_{m=1}^{\tilde M}.
\end{equation}
The Banach space $\M_{\Dd}(\mathrm{X})$ is defined as in \eqref{eq:locnative}. One can show that it admits a predual, as in Proposition \ref{prop:predual} . We now provide a list of useful assumptions.

\begin{itemize}
    \item []\hspace{-0.5cm}\emph{Assumption 1} The time points $(t_m)_{m=1}^M$ are such that $\forall m,\,t_m<t_{m+1}.$
    \item []\hspace{-0.5cm}\emph{Assumption 2} Sampling functionals are convex combinations of left and right limits, with the value of the parameter $\omega$ being in $[0,1].$
    \item []\hspace{-0.5cm}\emph{Assumption 3} \label{ass:4} The data-fidelity functional $E:\R^{M+\tilde M}\times\R^{M+\tilde M}\to\R^+\cup\{\infty\}$ is proper, coercive, lower-semicontinuous with its domain being an open set, and strictly convex in its second argument. In addition, it is lower-bounded, so that 
    $\underset{\mathbf{z}\in\R^{M+\tilde M}}{\text{inf}}E(\mathbf{y},\mathbf{z})=A>-\infty.
    $
    \item []\hspace{-0.5cm}\emph{Assumption 4} The measurement functional $\langle\bm{\nu},\cdot\rangle:
        \M_{\Dd}(\K)\to\R^{\tilde M}$ is $\w$-continuous.
    \item []\hspace{-0.5cm}\emph{Assumption 5} The underlying system $(1,\iota)$ is the unique ideal $\K$-localized system, as given in Theorem \ref{prop:ideal2}.
\end{itemize}
Our analysis of the solution set $\mathcal{V}(\K)$ is based on the description of a sequence of appropriate, approximate solution sets $\mathcal{V}(\K_{\epsilon})$. We define $\K_\epsilon=\K\setminus\K_\epsilon^c$ and 
\begin{equation}
    \K_\epsilon^c=\bigcup_{m=1}^{M}]t_m-\epsilon,t_m[\cup]t_m,t_m+\epsilon[.
\end{equation}
\begin{itemize}
    \item []\hspace{-0.5cm}\emph{Assumption 6} The parameter $\epsilon$ is small enough that $ t_{m-1}<(t_m-\epsilon)<(t_m+\epsilon)<t_{m+1}.$
    \item []\hspace{-0.5cm}\emph{Assumption 7} The closed set $\mathrm{supp}(\bm{\nu})=\bigcup_{m=1}^{\tilde M}\mathrm{supp}(\nu_m)$ is such that $\mathrm{supp}(\bm{\nu})\cap\{t_1,\ldots,t_M\}=\emptyset,$ and $\epsilon$ is small enough that $\mathrm{supp}(\bm{\nu})\cap\overline{\K_{\epsilon}^c}=\emptyset.$
\end{itemize}

\subsection{Special Case}
We establish the existence of solutions in the case where the supports of $\bb{t}$ and $\bm{\nu}$ do not overlap (Assumption 7). Let 
\begin{equation}
    \mathcal{J}^{\star}(\mathrm{X})=\underset{f\in\M_{\Dd}(\mathrm{X})}{\text{inf}}\mathcal{J}(f)
\end{equation}
be the optimal value of the loss functional. The proof of Theorem \ref{th:weakRT} is given in \ref{app:5.2}.
\begin{Theorem}
\label{th:weakRT}
If Assumptions 1 to 7 hold,
then the set $\mathcal{V}(\K_\epsilon)$ is nonempty, convex, and $\w$-compact. Its extreme points are of the form 
\begin{equation}
\label{eq:4.1.73}
    c+\sum_{k=1}^{K}c_ku(\cdot-x_k),\quad (\norm{c}_0+K)\leq(M+\tilde M),\quad x_k\in\K_{\epsilon}.
\end{equation}
In addition, the following holds:
\begin{equation}
\label{eq:4.1.74}
    \mathcal{J}^{\star}(\K_{\epsilon})=\mathcal{J}^{\star}(\K)\quad\quad\text{and}\quad\quad\bigcup_{\epsilon>0}\mathcal{V}(\K_{\epsilon})\subseteq\mathcal{V}(\K)\subseteq\overline{\mathcal{V}(\K)}^{\mathrm{weak}^{\star}}=\overline{\bigcup_{\epsilon>0}\mathcal{V}(\K_\epsilon)}^{\mathrm{weak}^{\star}},
\end{equation}
where the extreme points of  $\,\overline{\bigcup_{\epsilon>0}\mathcal{V}(\K_\epsilon)}^{\mathrm{weak}^{\star}}$ are as in \eqref{eq:4.1.73}, with $x_k\in\K$. Finally, if $\tilde{M}=0$, then one has the strict inclusion
\begin{equation}
    \mathcal{V}(\K)\subset\overline{\bigcup_{\epsilon>0}\mathcal{V}(\K_\epsilon)}^{\mathrm{weak}^{\star}}.
\end{equation}
\end{Theorem}

\begin{figure}[h!]
    \centering
    \includegraphics[width=0.3\linewidth]{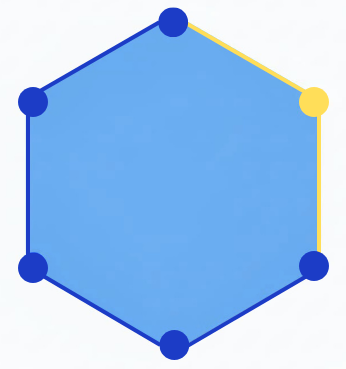}
    \caption{Schematic of the set $\overline{\mathcal{V}(\K)}^{\w}$. The shaded area corresponds to the interior $\overset{\circ}{\mathcal{V}(\K)}$ of the set, which is necessarily formed of minimizers. The darker (lighter, respectively) outline represents the part of the boundary that is formed by minimizers (non-minimizers, respectively).}
    \label{fig:3}
\end{figure}

Theorem \ref{th:weakRT} states that it is possible to extend the usual sparse extreme-point description of solution sets of variational inverse problems to the setting in which some of the measurement functionals are not $\w$-continuous but only continuous. In particular, this establishes that an optimization with measurements given by BV samples is a well-founded endeavor. 

The classic results \citep{bredies2020sparsity,flinth2019exact,unser2017splines} adapt as follows. The nonempty solution set $\mathcal{V}(\K)$ is not $\w$-closed (lighter dots in Figure \ref{fig:3}), but admits a $\Dd$-spline solution (darker dots in Figure \ref{fig:3}). Moreover, its closure (Figure \ref{fig:3}) possesses the same extreme-point representation as in the standard setting, where the measurement functional is $\w$-continuous.

\subsection{General Case}
In applications such as the Hawkes-process intensity estimation \citep{laub2015hawkes}, the loss functional involves the samples of a function as well as its integral on some intervals. We now relax Assumption 7 and treat the case where the supports of $\bb{t}$ and $\bm{\nu}$ overlap. We first illustrate that, in this general case, the existence of a solution is not guaranteed. 

\begin{Proposition}
\label{prop:counterexample}
The optimization problem as in \eqref{eq:3.1.30} with $\K=[0,1]$ and 
\begin{equation}
\label{eq:5.3.99}
    \mathcal{J}(f)=(\langle\delta_0^+,f\rangle)^2+(1-\langle\mathbbm{1}_{[0,1]},f\rangle)^2+\norm{f}_{\M_{\Dd}}
\end{equation}
has no solution.
\end{Proposition}

\begin{proof}[\textbf{Proof of Proposition \ref{prop:counterexample}}]
Observe that \eqref{eq:5.3.99} leads to a solution if and only if 
\begin{equation}
\label{eq:5.3.100}
    \underset{f\in\Dd^{-1}\{\M([0,1])\}}{\mathrm{argmin}}\,(\langle\delta_0^+,f\rangle)^2+(1-\langle\mathbbm{1}_{[0,1]},f\rangle)^2+\norm{\Dd\{f\}}_{\M}
\end{equation}
has a solution. Assume, by contradiction, that \eqref{eq:5.3.100} has a solution $f^{\star}=u\ast m^{\star}$. If $m^{\star}$ is nonzero over $]\epsilon,1]$, with $\epsilon>0$, then 
\begin{equation}
    f^{\star\star}=u\ast(m^{\star}\mathbbm{1}_{[0,\epsilon]}+\delta_{\epsilon}a)\quad\text{with}\quad 
    a=\frac{1}{(1-\epsilon)}\int_{\epsilon}^1m^{\star}(]\epsilon,x])\,\mathrm{d}x
\end{equation}
is such that 
\begin{equation}
    (\langle\delta_0^+,f^{\star}\rangle)^2+(1-\langle\mathbbm{1}_{[0,1]},f^{\star}\rangle)^2=(\langle\delta_0^+,f^{\star\star}\rangle)^2+(1-\langle\mathbbm{1}_{[0,1]},f^{\star\star}\rangle)^2.
\end{equation}
Then, 
\begin{equation}
    \Big\vert\int_{\epsilon}^1m^{\star}(]\epsilon,x])\,\mathrm{d}x\Big\vert<(1-\epsilon)\vert m^{\star}(]\epsilon,1])\vert\Rightarrow\norm{\Dd\{f^{\star\star}\}}_{\M}<\norm{\Dd\{f^{\star}\}}_{\M}.
\end{equation}
Therefore, $\forall\epsilon>0$, one must have that $\mathrm{supp}(m^{\star})\subset[0,\epsilon]$ and, consequently, that $\mathrm{supp}(m^{\star})=\{0\}.$ Thus,
\begin{equation}
    f^{\star}(\cdot)=\alpha u(\cdot)\quad\text{with}\quad\mathcal{J}(f^{\star})=\alpha^2+(1-\alpha)^2+\alpha.
\end{equation}
Finally, we consider $\tilde f=\frac{\alpha}{(1-\epsilon)} u(\cdot-\epsilon)$ and calculate that 
\begin{equation}
    \mathcal{J}(\tilde f)=0+(1-\alpha)^2+\frac{\alpha}{(1-\epsilon)},
\end{equation}
which is smaller than $\mathcal{J}(f^{\star})$ for $\epsilon$ sufficiently small. This contradicts the optimality of $f^{\star}.$
\end{proof}

Although the solution set $\mathcal{V}(\K)$ may be empty, the optimal value $\mathcal{J}^{\star}(\K)$ can always be approximated by a sequence of $\Dd$-splines that are solutions of the optimization problems posed over $\K_{\epsilon}$. The proof of Theorem \ref{th:weakRT2} is given in \ref{app:5.3}.

\begin{Theorem}
\label{th:weakRT2}
Under Assumptions 1 to 6, the set $\mathcal{V}(\K_\epsilon)$ is nonempty, convex, and $\w$-compact. Its extreme points are of the form 
\begin{equation}
    c+\sum_{k=1}^{K}c_ku(\cdot-x_k),\quad (\norm{c}_0+K)\leq M,\quad x_k\in\K_{\epsilon}.
\end{equation}
In addition, the following convergence holds:
\begin{equation}
\label{eq:4.2.125}
    \lim_{\epsilon\to0}\mathcal{J}^{\star}(\K_{\epsilon})=\mathcal{J}^{\star}(\K).
\end{equation}
\end{Theorem}
More precisely, the net $(f_{\epsilon}^{\star})_{\epsilon}$, with $f_{\epsilon}^{\star}\in\mathcal{V}(\K_{\epsilon})$ and that approximates $\mathcal{J}^{\star}(\K)$, can be chosen  $\w$-convergent to some limit $\tilde f$. In practice, this can be sufficient for the practitioner as the sequence $(f_{\epsilon}^{\star})_{\epsilon}$ is convergent and approximates the optimal value $\mathcal{J}^{\star}(\K)$ arbitrarily well, its sole constraint being to have no knots closer than $\epsilon$ to $\bb{t}$. In fact, this constraint can easily be implemented by discretizing the problem on grids that do not overlap with $\bb{t}$ \citep{guillemet2025convergence}.

Moreover, one can relax the condition on $(f_{\epsilon}^{\star})_{\epsilon}$ by requiring only the convergence of a subsequence. The set of such limits $\tilde f$ is called the Kuratowski limit superior $\bb{Ls}(\mathcal{V}(\K_{\epsilon}))$, defined as \citep{kuratowski2014topology}
\begin{equation}
    \bb{Ls}(\mathcal{V}(\K_{\epsilon}))=\big\{\tilde f\in\M_{\Dd}(\K):\quad \tilde f=\w-\underset{\ell\to\infty}{\mathrm{lim}}f_{\epsilon_{\ell}}^{\star},\quad f_{\epsilon_{\ell}}^{\star}\in\mathcal{V}(\K_{\epsilon_{\ell}})\big\}.
\end{equation}

Kuratowski’s notions of limit superior and limit inferior are closely related to the concept of $\Gamma$-convergence \citep[Definition~4.1]{dal2012introduction}, which is itself used to characterize the convergence of solution sets in regularized inverse problems \citep{duval2017sparse,fageot2025variational}. However, this stronger tool cannot be applied in our setting, as some of the basic assumptions required for $\Gamma$-convergence can not be satisfied. Consequently, we rely solely on the limit superior $\bb{Ls}(\mathcal{V}(\K_{\epsilon}))$, whose characterization is provided in Theorem~\ref{th:kuratowski}. The proof is given in \ref{app:5.3}.

\begin{Theorem}
\label{th:kuratowski}
 If Assumptions 1 to 6 hold, then $\bb{Ls}(\mathcal{V}(\K_{\epsilon}))$ is nonempty, and $\w$-compact. Its extreme points are all of the form 
 \begin{equation}
 \label{eq:5.3.107}
    c+\sum_{k=1}^{K}c_ku(\cdot-x_k),\quad (\norm{c}_0+K)\leq M,\quad x_k\in\K.
\end{equation}
Finally, for every $\tilde f\in\bb{Ls}(\mathcal{V}(\K_{\epsilon}))$, there exists a net $(f_{\epsilon}^{\star})_{\epsilon}$ such that $\underset{\epsilon\to0}{\mathrm{lim}}\mathcal{J}(f^{\star}_{\epsilon})=\mathcal{J}^{\star}(\K)$, which admits a sequence $(f_{\epsilon_{\ell}}^{\star})_{\ell=1}^{\infty}$ that is $\w$-convergent to $\tilde f$.
\end{Theorem}
To better visualize the relationship between  $\bb{Ls}(\mathcal{V}(\K_{\epsilon}))$ and $\mathcal{V}(\K)$, we provide an itemized comparison.

\begin{itemize}
    \item [1.] If the measurement operator is $\w$-continuous, then $\mathcal{V}(\K)=\bb{Ls}(\mathcal{V}(\K_{\epsilon}))$, and both are nonempty, convex, and $\w$-compact, with extreme points of the form \eqref{eq:5.3.107}.
    \item [2.] If the measurement operator is the concatenation of the sampling functional $\delta_{t}^{\omega}$ with $\w$-continuous functionals whose supports do not overlap, then $\mathcal{V}(\K)\subset\bb{Ls}(\mathcal{V}(\K_{\epsilon}))$. The solution set is nonempty and convex but may not be $\w$-compact. The limit superior is nonempty, convex, and $\w$-compact, with extreme points of the form \eqref{eq:5.3.107}.
    \item [3.] If the measurement operator is the concatenation of the sampling functional $\delta_{t}^{\omega}$ with $\w$-continuous functionals whose supports do overlap, then $\mathcal{V}(\K)\subset\bb{Ls}(\mathcal{V}(\K_{\epsilon}))$. The solution set may be empty. The limit superior is nonempty, and $\w$-compact, with extreme points of the form \eqref{eq:5.3.107}.
\end{itemize}

In Item 1, the equality $\mathcal{V}(\K)=\bb{Ls}(\mathcal{V}(\K_{\epsilon}))$ follows from the fact that, in this case, $\K=\K_{\epsilon}$. In Item 2, the embedding $\mathcal{V}(\K)\subset\bb{Ls}(\mathcal{V}(\K_{\epsilon}))$ follows from the equality 
\begin{equation}
\bb{Ls}(\mathcal{V}(\K_{\epsilon}))=\overline{\bigcup_{\epsilon>0}\mathcal{V}(\K_\epsilon)}^{\text{weak}^{\star}},
\end{equation}
whose proof is straightforward. 

The fewer properties are assumed of the measurement functional, the more properties are lost in the solution set. Nevertheless, these properties are retained by the limit superior set of solutions, which shows that they generalize the solution set and make it an appropriate replacement. Importantly, $\bb{Ls}(\mathcal{V}(\K_{\epsilon}))$ corresponds to the set of all limits to which one can converge while attempting to solve the optimization problem iteratively on $\M_{\Dd}(\K_\epsilon)$.

We illustrate Theorem~\ref{th:kuratowski} by extending the example presented in Proposition~\ref{prop:counterexample}. For $\epsilon < 0.5$, the solution in $\mathcal{V}(\K_{\epsilon}) = \mathcal{V}([\epsilon,1])$ is unique and given by
\begin{equation}
    \mathcal{V}(\K_{\epsilon})
    = \left\{ \frac{1 - 2\epsilon}{2(1-\epsilon)^2}\, u(\cdot - \epsilon) \right\},
    \qquad
    \text{with}
    \qquad
    \mathcal{J}^{\star}(\K_{\epsilon})
    = \frac{3 - 4\epsilon}{4(1-\epsilon)^2}.
\end{equation}

\noindent Consequently, the solution in $\bb{Ls}(\mathcal{V}(\K_{\epsilon}))$ is also unique and satisfies
\begin{equation}
    \bb{Ls}(\mathcal{V}(\K_{\epsilon}))
    = \left\{ \tfrac{1}{2}\, u(\cdot) \right\},
    \qquad
    \text{with}
    \qquad
    \mathcal{J}\!\left( \tfrac{1}{2} u(\cdot) \right)
    = 1
    <
    \mathcal{J}^{\star}
    = \tfrac{3}{4}.
\end{equation}


\section{Higher-Order BV Spaces}
With the same kind of techniques, we can extend the results of Sections 4 and 5 to the Banach space of $N^{\mathrm{th}}$-order BV functions. For concision, we list the results without giving proofs.

\subsection{Banach Spaces}
Let $\Dd^N$ be the $N^{\mathrm{th}}$-derivative operator, with unique causal Green's function 
$g_N(\cdot)=\frac{(\cdot)_+^{N-1}}{(N-1)!}$ and null space 
$\mathcal{N}_{N}=\mathrm{span}\left\{\frac{(\cdot)^{n-1}}{(n-1)!}\right\}_{n=1}^N$ of dimension $N$. The space of functions whose $N^{\mathrm{th}}$ weak derivative is a bounded measure is defined as 
\begin{equation}
    \M_{N}(\R) = \bigl\{ f:\R\to\R \quad\text{s.t.}\quad\norm{\Dd^N\{f\}}_{\M}<\infty\bigr\},
\end{equation}
with $\M_{1}(\R)$ corresponding to $\M_{\Dd}(\R)$. As in the BV case, the specification of an inverse $\Dd^{-N}$ requires the introduction of a bi-orthogonal system $(\bb{p}, \bm{\phi})$.

\begin{Definition}
Let $\bb{p}=(p_n)_{n=1}^N$ be a basis of $\mathcal{N}_{N}$ and $\bm{\phi}=(\phi_{n})_{n=1}^N$ be a complementary set of analysis functions. Then, $(\bb{p}, \bm{\phi})$ is a universal $\DN$-admissible system if, for all $n\in[1\ldots N]$, the function $\phi_n\in\mathcal{D}(\R)$ is compactly supported with $\mathrm{supp}(\phi_n)\subset[\phi^-,\phi^+]$ and satisfies 
\begin{equation}
    \forall (n,m)\in\{1,\ldots,N\}^2:\quad\langle 
    p_n, \phi_m\rangle=\delta_{n-m}.
\end{equation}
\end{Definition}

The associated kernel $g_{\bm{\phi}}$ and integral transform $\Dd_{\bm{\phi}}^{-N}$ are such that \citep{unser2017splines}
\begin{equation}
\label{eq:2.8.2}
  g_{\bm{\phi}}(t,x)=g_N(t-x)-\sum_{n=1}^N\langle\phi_{n}(\cdot), g_N(\cdot-x)\rangle p_n(t), 
\end{equation}
\begin{equation}
\label{eq:2.9}
     \Dd_{\bm{\phi}}^{-N}\{f\}=\int_{\R}g_{\bm{\phi}}(\cdot,x)f(x)\,\mathrm{d}x.
\end{equation}
This integral transform yields the direct-sum decomposition 
\begin{equation}
    \M_{N}(\R)=\Dd_{\bm{\phi}}^{-N}\{\M(\R)\}\oplus\mathcal{N}_{N},
\end{equation}
where $\M_{N}(\R)$ is a Banach space for the norm
\begin{equation}
    \norm{f}_{\M_{N}}:=\norm{\DN\{f\}}_{\M}+\sum_{n=1}^N\vert\langle\phi_n,f\rangle\vert.
\end{equation}

Observe that $\bm{\phi}$ forms a basis for the dual space $\mathcal{N}_{N}^{\star}$ of $\mathcal{N}_{N}$. The Banach space 
$\C_{N}(\R)=\Dd^{N,\star}\{\C_0(\R)\}\oplus\mathcal{N}_N^{\star}$ is a predual for $\M_{N}(\R)$ \citep{unser2017splines}. This characterization enables us to extend the notion of an admissible system to a collection $(\phi_n)_{n=1}^N$ with $\phi_{n}\in\C_{N}(\R)$ \citep{unser2019native}.

The definition of $\M_N(\R)$ motivates the introduction of generalized sampling functionals.

\subsection{Generalized Sampling}

We recall that the adjoint of $\Dd^n$ is $\Dd^{n\star}=(-1)^n\Dd^n$. We split generalized sampling functionals into two categories.
\begin{itemize}
    \item The first category consists of functionals spanned by $\Dd^{n}\{\delta_t\}$, with $n<(N-1)$ and $t\in\R$.
    \item The second category consists of functionals spanned by $\Dd^{N-1}\{\delta_t\}$, with $t\in\R$.
\end{itemize}
Functionals in the first category are the more regular ones. Their continuity is studied in Proposition \ref{diractweakstar}, whose proof follows directly from the observation that 
$\Dd^{n}\{\delta_{t}\}\in\C_{\Ls}(\R)$ \citep[Corollary 1]{guillemet2025adaptive}.  

\begin{Proposition}
\label{diractweakstar}
For all $n\leq(N-2)$, the mapping 
\begin{equation}
\langle\Dd^{n\star}\{\delta_{t}\},\cdot\rangle:\begin{cases}
    \M_{N}(\R)\to\R,\\
    f\mapsto\langle\Dd^{n\star}\{\delta_{t}\},f\rangle=\langle\delta_t,\Dd^{n}\{f\}\rangle=f^{(n)}(t)
\end{cases}       
\end{equation}
is well-defined, linear, and $\w$-continuous. 
\end{Proposition}

The study of functionals in the second category is more delicate, since $f^{(N-1)}$ is only well-defined via its left and right traces $(f^{(N-1),-},f^{(N-1),+})$. We introduce the generalized left- and right-limit functionals, $\forall f\in\M_{\Dd}(\R)$,
\begin{equation}
    \langle\Dd^{N-1\star}\{\delta_{t}^{-}\},f\rangle=\underset{\epsilon\downarrow0}{\mathrm{lim }}f^{(N-1)}(t-\epsilon),\quad\text{and}\quad\langle\Dd^{N-1\star}\{\delta_{t}^{+}\},f\rangle=\underset{\epsilon\downarrow0}{\mathrm{lim }}f^{(N-1)}(t+\epsilon),
\end{equation}
and further define the convex combination, for $\omega\in[0,1]$,
\begin{equation}
  \Dd^{N-1\star}\{\delta_{t}^{\omega}\}=\omega\Dd^{N-1\star}\{\delta_{t}^{+}\}+(1-\omega)\Dd^{N-1\star}\{\delta_{t}^{-}\}.
\end{equation}
The continuity of $\Dd^{N-1\star}\{\delta_{t}^{\omega}\}$ is discussed in Proposition \ref{diractweak3}.
\begin{Proposition}
\label{diractweak3}
The mapping 
\begin{equation}
\langle\Dd^{N-1\star}\{\delta_{t}^{\omega}\},\cdot\rangle:\begin{cases}
    \M_{N}(\R)\to\R,\\
    f\mapsto\langle\Dd^{N-1\star}\{\delta_{t}^{\omega}\},f\rangle
\end{cases}       
\end{equation}
is well-defined, linear, continuous, but is not $\w$-continuous.
\end{Proposition}

\subsection{Fundamental Systems}

Based on our understanding of the continuity of $\Dd^{n\star}\{\delta_{t}^{\omega}\}$, we proceed as in Section 4 to define a unique local inverse $\Dd^{-N}$ of $\Dd^{N}.$ We start by extending the notion of localized systems.

\begin{Definition}
    A $\DN$-admissible system $(\bb{p},\bm{\phi})$ is called $\K$-localized if $[\phi^-,\phi^+]\subset[K^-,K^+]$.
\end{Definition}

The $\K$-localized systems yield a favorable weight distribution of the kernel $g_{\bm{\phi}}$. Indeed, an analogous result to Proposition \ref{prop:localizedjustification} can be established. Such systems lead to the definition of $\M_{N}(\K)$ as \citep{guillemet2025convergence}
\begin{align}
    \M_{N}(\K)&=\{f\in\M_{N}(\R):\quad\Dd^{N}\{f\}\in\M(\K)\}\nonumber\\
    &=\Dd^{-N}_{\bm{\phi}}\{\M(\K)\}\oplus\mathcal{N}_{N},
\end{align}
in which every function $f$ admits a non-direct-sum representation (see Section 4.1) of the form $f=(g_N\ast\mu)+p$ with $\mu\in\M(\K)$ and $p\in\mathcal{N}_N.$ 

To make this sum direct, and to construct a kernel $g_{\bm{\phi}}$ that is locally shift-invariant, we extend the previously introduced notion of fundamental systems. To this end, it is useful to remark that the same construction as in Proposition \ref{prop:predual} yields that $\M_{N}(\K)$ has a predual.

\begin{Definition}
    A $\K$-localized system $(\bb{p},\bm{\phi})$ is called $\K$-fundamental if
    \begin{itemize}
        \item[1.] for all $n\in[1\ldots N]$, the functional $\phi_n:\M_{N}(\K)\to\R$ is $\w$-continuous;
        \item[2.] the kernel $g_{\bm{\phi}}$ reduces to the Green’s function, in the sense that
    \begin{equation}
    \forall x\in[K^-,\infty[:\quad g_{\bm{\phi}}(\cdot,x)=g_{N}(\cdot-x);
    \end{equation}
        \item[3.] the null-space element $p_{n}$ is the shifted canonical polynomial
    \begin{equation}
    \label{eq:4.4.71}
        p_{n}(\cdot)=\frac{(\cdot-K^-)^{n-1}}{(n-1)!}.
    \end{equation}
    \end{itemize}
\end{Definition}

To tell it apart from $(\bb{p},\bm{\phi})$, we denote a $\K$-fundamental system by $(\bb{p},\bm{\iota})$. As before, this choice of kernel leads to a significant simplification of the theory, with
\begin{equation}
    \M_{N}(\K)=\Dd_{\bm{\iota}}^{-N}\{\M(\K)\}\oplus\mathcal{N}_{N}=\Dd^{-N}\{\M(\K)\}\oplus\mathcal{N}_{N},
\end{equation}
where the integral operator $\Dd^{-N}$ is described by the shift-invariant kernel $g_{N}$. 
As in the first-order case, we can show that there exists a unique $\K$-fundamental system.

\begin{Theorem}
 \label{prop:ideal3}
The system $(\bb{p},\bm{\iota})$ defined as follows is the unique $\K$-fundamental system: the null-space element $p_{i,n}$ is given by \eqref{eq:4.4.71}, and $\iota_{i,n}$ is defined by
\begin{align}
   \forall n\in[1\ldots N-1],\quad &\iota_{n}=\Dd^{n-1\star}\{\delta_{K^-}\},\quad\text{and}\quad\iota_{N}=\Dd^{N-1\star}\{\delta_{K^-}^-\}.
\end{align}
\end{Theorem}

Theorem \ref{prop:ideal3} shows that there exists a unique, natural norm over $\M_{N}([0,T])$, which is given by the standard 
\begin{equation}
    \norm{f}_{\M_N}=\norm{\Dd^N\{f\}}_{\M}+\vert f^{(N-1)}(0^-)\vert+\sum_{n=1}^{N-1}\vert f^{(n-1)}(0)\vert.
\end{equation}

\subsection{Resolution of Optimization Problems}
We are now interested in solution sets of the form
\begin{equation}
\label{eq:3.1.30'}
\mathcal{J}^{\star}(\mathrm{X})=\underset{f\in\M_{N}(\mathrm{X})}{\text{inf}}\mathcal{J}(f),\quad\quad\mathcal{V}(\mathrm{X})=\underset{f\in\M_{N}(\mathrm{X})}{\text{argmin}}\quad\mathcal{J}(f),
\end{equation}
with the \emph{loss functional}
\begin{equation}
\mathcal{J}(f)=E\left(\bb{y},\left\langle\begin{bmatrix}
\Dd^{N-1,\star}\{\bm{\bm{\delta}_{\bb{t}}}^{\omega}\}\\
\bm{\nu}
\end{bmatrix},f\right\rangle\right)+\lambda\norm{f}_{\M_{N}}.
\end{equation}

The new Assumptions 1', 2', 3', 5', and 6' are identical to Assumptions 1, 2, 3, 5, and 6. Assumption 4' is formulated as follows.
\begin{itemize}
    \item []\hspace{-0.5cm}\emph{Assumption 4'} The measurement functional $\langle\bm{\nu},\cdot\rangle:
        \M_{N}(\K)\to\R^{\tilde M}$ is $\w$-continuous.
\end{itemize} 

In Theorem \ref{th:weakRT3}, we investigate the existence of solutions as well as the structure of the solution set. The division into cases based on overlapping and non-overlapping support of $\bb{t}$ and $\mathrm{supp}(\bm{\nu})$ is no longer appropriate. Instead, we consider the more fundamental distinction between $\tilde M=0$ and $\tilde M\neq0$.

\begin{Theorem}
\label{th:weakRT3}
Under Assumptions 1' to 6',
the set $\mathcal{V}(\K_\epsilon)$ is nonempty, convex and $\w$-compact. Its extreme points are of the form 
\begin{equation}
\label{eq:4.1.73.2}
    \sum_{n=1}^Nc_n\frac{(\cdot)^{n-1}}{(n-1)!}+\sum_{k=1}^{K}a_k\frac{(\cdot-x_k)_+^{N-1}}{(N-1)!},\quad (\norm{\bb{c}}_0+K)\leq(M+\tilde M),\quad x_k\in\K_{\epsilon}.
\end{equation}
In addition, 
\begin{itemize}
    \item if $\tilde M=0$, then 
    \begin{equation}
    \mathcal{J}^{\star}(\K_{\epsilon})=\mathcal{J}^{\star}(\K)
\end{equation}
and 
\begin{equation}
    \bigcup_{\epsilon>0}\mathcal{V}(\K_{\epsilon})\subseteq\mathcal{V}(\K)\subset\overline{\bigcup_{\epsilon>0}\mathcal{V}(\K_\epsilon)}^{\mathrm{weak}^{\star}}=\bb{Ls}(\mathcal{V}(\K_{\epsilon}));
\end{equation}
    \item if $\tilde M\neq0$, then for every $\tilde f\in\bb{Ls}(\mathcal{V}(\K_{\epsilon}))$ there exists a net $(f_{\epsilon}^{\star})_{\epsilon}$ with 
    $f_{\epsilon}^{\star}\in\mathcal{V}(\K_{\epsilon})$ such that 
    \begin{equation}
        \underset{\epsilon\to0}{\mathrm{lim}}\,\mathcal{J}(f^{\star}_{\epsilon})=\mathcal{J}^{\star}(\K)\Rightarrow     \lim_{\epsilon\to0}\mathcal{J}^{\star}(\K_{\epsilon})=\mathcal{J}^{\star}(\K), 
    \end{equation}
    and such that it admits a sequence $(f_{\epsilon_{\ell}}^{\star})_{\ell=1}^{\infty}$ that is $\w$-convergent to $\tilde f$.
\end{itemize}
In addition, in both cases the set $\bb{Ls}(\mathcal{V}(\K_{\epsilon}))$ is nonempty, $\w$-compact, with extreme points of the form \eqref{eq:4.1.73.2}, with $x_k\in\K$.
\end{Theorem}

Theorem \ref{th:weakRT3} generalizes Theorem \ref{th:weakRT}, Theorem \ref{th:weakRT2}, and Theorem \ref{th:kuratowski} simultaneously. The results are analogous to the $BV$ case. 

When $\tilde{M}=0$, Theorem \ref{th:weakRT3} establishes the well-posedness of variational inverse problems involving generalized $BV$ samples. In particular, it shows that there always exists a $\DN$-spline solution. In this case, the solution set is not $\w$-closed, but its closure $\bb{Ls}(\mathcal{V}(\K_{\epsilon}))$ has the same extreme-point structure as the solution set in the $\w$-continuous measurement scenario \citep{unser2017splines}.

When $\tilde{M}\neq0$, the solution set may be empty, but one can always approximate the optimal value $\mathcal{J}^{\star}(\K)$ with a $\w$-convergent sequence of $\DN$-splines. The set of limits of such sequences is $\bb{Ls}(\mathcal{V}(\K_{\epsilon}))$, which admits, as before, the same extreme-point structure as the solution set in the $\w$-continuous measurement scenario.

\section*{Acknowledgment}
Vincent Guillemet was supported by the Swiss National Science Foundation (SNSF) under Grant 200020\_219356.

\section*{Data availability}
Data sharing is not applicable to this article as no datasets were generated or analysed.
\appendix
\section{Proofs}
\subsection{Proofs of Section 5.2}
\label{app:5.2}

\begin{proof}[\textbf{Proof of Theorem \ref{th:weakRT}}]
\hfill\\
\textbf{Step 1.} We study $\mathcal{V}(\K_\epsilon)$. The measurement functional  
\begin{equation}
\langle\bm{\delta}_{\bb{t}}^\omega,\cdot\rangle:
    \M_{\Dd}(\K_{\epsilon})\to\R^M
\end{equation}
is $\w$-continuous for the same reason as in Proposition \ref{diracweakstar2}. By the direct method of the calculus of variations and the Krein–Milman theorem, we obtain that $\mathcal{V}(\K_\epsilon)$ is nonempty, convex, $\w$-compact, and is the $\w$-closed convex hull of its extreme points. Finally, \citep[Corollary 3.7]{boyer2019representer} shows that the extreme points of $\mathcal{V}(\K_\epsilon)$ are of the desired form.
\hfill\\
\textbf{Step 2.}
We prove that $\mathcal{V}(\K_\epsilon)\subset\mathcal{V}(\K)$. For $f=c+u\ast \mu\in\M_{\Dd}(\K),$ we define the projection operator 
\begin{align}
    \mathrm{proj}_{\K_{\epsilon}}(\cdot)\begin{cases}
    \M_{\Dd}(\K)\mapsto\M_{\Dd}(\K_{\epsilon})\\
    f\mapsto c+(u\ast \mu_{\epsilon}),
    \end{cases}
\end{align}
where
\begin{equation}
    \mu_{\epsilon}=\mu\mathbbm{1}_{\K_\epsilon}+\sum_{m=1}^M\Big(a_m\delta_{t_m-\epsilon}+b_m\delta_{t_m+\epsilon}\Big),\quad\text{with}\quad\begin{cases}
            a_m=\mu(]t_m-\epsilon,t_m[),\\
            b_m=\mu(]t_m,t_m+\epsilon[).
    \end{cases}
\end{equation}
For every $f\in\M_{\Dd}(\K)$, our claims are the following:
\begin{itemize}
    \item [1.] the inequality $\norm{f_{\epsilon}}_{\M_{\Dd}}\leq\norm{f}_{\M_{\Dd}}$ holds; 
    \item [2.] for all $t\in\bigcup_{m=1}^{M}\{t_m\}\cup]t_{m}+\epsilon,t_{m+1}-\epsilon[$,
    \begin{equation}
    \left\langle\delta_{t}^-,u\ast\mu_{\epsilon}\right\rangle=\langle\delta_{t}^-,u\ast \mu\rangle,\quad
    \text{and}\quad\left\langle\delta_{t}^+,u\ast\mu_{\epsilon}\right\rangle=\langle\delta_{t}^+,u\ast \mu\rangle.
    \end{equation}
\end{itemize}
We only prove Claim 2, since Claim 1 is straightforward. We show Claim 2 for $t\in\{t_m\}_{m=1}^M$ and for right limits. From Proposition \ref{prop:weakrep}, it follows that 
\begin{align}
\left\langle\delta_{t_m}^+,u\ast\mu_{\epsilon}\right\rangle=&\mu([K^-,t_m]\cap\K_{\epsilon})\nonumber\\
&+\sum_{m'=1}^{M}\Big((t_m-t_{m'}+\epsilon)_{+}^0a_{m'}+(t_m-t_{m'}-\epsilon)_{+}^0b_{m'}\Big)\nonumber\\
=&\mu([K^-,t_m]\cap\K_{\epsilon})+\mu(]t_{m}-\epsilon,t_{m}[)\nonumber\\
&+\sum_{m'=1}^{m-1}\Big(\mu(]t_{m'}-\epsilon,t_{m'}[)+\mu(]t_{m'},t_{m'}+\epsilon[)\Big)\nonumber\\
=&\mu([K^-,t_m])=\langle\delta_{t_m}^+,u\ast \mu\rangle.
\end{align}
Claim 2 therefore holds, and implies that 
\begin{align}
     \left\langle\begin{bmatrix}
\bm{\bm{\delta}_{\bb{t}}}^{\omega}\\
\bm{\nu}
\end{bmatrix},f\right\rangle=\left\langle\begin{bmatrix}
\bm{\bm{\delta}_{\bb{t}}}^{\omega}\\
\bm{\nu}
\end{bmatrix},\mathrm{proj}_{\K_{\epsilon}}(f)\right\rangle&\Rightarrow \mathcal{J}(\mathrm{proj}_{\K_{\epsilon}}(f))\leq\mathcal{J}(f)\nonumber\\
&\Rightarrow\mathcal{J}^{\star}(\K_{\epsilon})\leq\mathcal{J}^{\star}(\K)\nonumber\\
&\Rightarrow\mathcal{V}(\K_{\epsilon})\subset\mathcal{V}(\K),
\end{align}
where we used the embedding
 $\M_{\Dd}(\K_{\epsilon})\subset\M_{\Dd}(\K)$.
\hfill\\
 \textbf{Step 3.}
We prove that $\mathcal{V}(\K)\subset\overline{\bigcup_{\epsilon>0}\mathcal{V}(\K_\epsilon)}^{\text{weak}^{\star}}.$ To this end, for any $f^{\star}\in\mathcal{V}(\K)$, we construct a net $(f_{\epsilon}^{\star})_{\epsilon}$ that is $\text{weak}^{\star}$-convergent to $f^{\star}$ and satisfies $f_{\epsilon}^{\star}\in\mathcal{V}(\K_{\epsilon})$.

Define $f_{\epsilon}^{\star}=\mathrm{proj}_{\K_{\epsilon}}(f)$ and observe that, by the argument of Step 2, $f_{\epsilon}^{\star}\in\mathcal{V}(\K_{\epsilon}).$ We now note that $(f_{\epsilon}^{\star})_{\epsilon}$ is $\text{weak}^{\star}$-convergent to $f^{\star}$ if and only if, for all $v\in\C_0(\R),$
\begin{align}
    \lim_{\epsilon\to0}\,\langle u\ast \mu_{\epsilon},\Dd^{\star}\{v\}\rangle=\langle u\ast \mu,\Dd^{\star}\{v\}\rangle\quad\Leftrightarrow\quad\lim_{\epsilon\to0}\,\langle \mu_{\epsilon},v\rangle=\langle \mu,v\rangle.
\end{align}
It is straightforward to verify that $\mu_{\epsilon}$ is indeed $\w$-convergent to $\mu$. 
\hfill\\
\textbf{Step 4.} We prove that the extreme points of $\overline{\bigcup_{\epsilon>0}\mathcal{V}(\K_\epsilon)}^{\text{weak}^{\star}}$ are of the desired form.
First, we observe that 
\begin{align}
&\mathcal{V}(\K_\epsilon)=\overline{\mathrm{co}\Big(\mathrm{ext}(\mathcal{V}(\K_{\epsilon}))\Big)}^{\w}\label{eq:4.1.90}\\
\Rightarrow&\mathcal{V}(\K_\epsilon)\subset\overline{\mathrm{co}\Big(\bigcup_\epsilon\mathrm{ext}(\mathcal{V}(\K_{\epsilon}))\Big)}^{\w}\label{eq:4.1.91}\\
\Rightarrow&\overline{\bigcup_{\epsilon}\mathcal{V}(\K_\epsilon)}^{\text{weak}^{\star}}\subset\overline{\mathrm{co}\Big(\bigcup_\epsilon\mathrm{ext}(\mathcal{V}(\K_{\epsilon}))\Big)}^{\w}\label{eq:4.1.92},
\end{align}
where \eqref{eq:4.1.90} follows from the Krein–Milman theorem, and \eqref{eq:4.1.92} follows from \eqref{eq:4.1.91} because the RHS of \eqref{eq:4.1.91} is $\w$-closed.
Second, we find that 
\begin{align}
    &\bigcup_{\epsilon}\mathrm{ext}(\mathcal{V}(\K_{\epsilon}))\subset\bigcup_{\epsilon}\mathcal{V}(\K_{\epsilon})\subset\overline{\bigcup_{\epsilon}\mathcal{V}(\K_\epsilon)}^{\text{weak}^{\star}}\label{eq:4.1.93}\\
    \Rightarrow\,&\overline{\mathrm{co}\Big(\bigcup_\epsilon\mathrm{ext}(\mathcal{V}(\K_{\epsilon}))\Big)}^{\w}\subset\overline{\bigcup_{\epsilon}\mathcal{V}(\K_\epsilon)}^{\text{weak}^{\star}}\label{eq:4.1.94},
\end{align}
where \eqref{eq:4.1.94} follows from \eqref{eq:4.1.93} because its RHS is convex and $\w$-closed. This shows that
\begin{equation}
\label{eq:claim3}
\overline{\mathrm{co}\Big(\bigcup_\epsilon\mathrm{ext}(\mathcal{V}(\K_{\epsilon}))\Big)}^{\w}=\overline{\bigcup_{\epsilon}\mathcal{V}(\K_\epsilon)}^{\text{weak}^{\star}}.
\end{equation}
Finally, together with \eqref{eq:claim3}, Milman's (partial) converse to the Krein–Milman theorem \citep{milman1947} implies that, because 
\begin{equation}
    \bigcup_\epsilon\mathrm{ext}(\mathcal{V}(\K_{\epsilon}))\subset\overline{\bigcup_{\epsilon}\mathcal{V}(\K_\epsilon)}^{\text{weak}^{\star}},
\end{equation}
whose RHS is $\w$-compact, one has that 
\begin{equation}
    \label{eq:4.1.96}
    \mathrm{ext}\Big(\overline{\bigcup_{\epsilon}\mathcal{V}(\K_\epsilon)}^{\text{weak}^{\star}}\Big)\subset\overline{\bigcup_\epsilon\mathrm{ext}(\mathcal{V}(\K_{\epsilon}))}^{\w}.
\end{equation}
The spline representation of the extreme points of $\overline{\bigcup_{\epsilon}\mathcal{V}(\K_\epsilon)}^{\text{weak}^{\star}}$ follows directly from \eqref{eq:4.1.96} and the representation of the extreme points of $\mathcal{V}_{\epsilon}$. 
\hfill\\
\textbf{Step 5.} We prove that the inclusion of Step 3 is strict in the case of left limits $(\omega=0)$. Let $t_0=K^-.$ If $f^{\star}=c+(u\ast\mu)\in\mathcal{V}(\K)$, then 
\begin{equation}
    f^{\star}_{\epsilon}=c+\sum_{m=0}^{M-1}\mu([t_m,t_{m+1}[)u(\cdot-t_{m+1}+\epsilon)
\end{equation}
is such that, for all $\epsilon\geq0$,
\begin{itemize}
    \item [1.] the inequality $\norm{f^{\star}_{\epsilon}}_{\M_{\Dd}}\leq\norm{f^{\star}}_{\M_{\Dd}}$ holds;
    \item [2.] for all $\epsilon>0$, $\langle\delta_{t_{m+1}}^-,f^{\star}_\epsilon\rangle=\langle\delta_{t_{m+1}}^-,f^{\star}\rangle=\mu([K^-,t_{m+1}[)$.
\end{itemize}
Therefore, $f_{\epsilon}^\star\in\mathcal{V}(\K_{\epsilon})$. Furthermore, it is straightforward to show that $f_{\epsilon}^\star$ is $\w$-convergent to $f_0^{\star}$. It remains to prove that $f_0^{\star}\notin\mathcal{V}(\K).$ To this end, observe that 
\begin{equation}
\label{eq:4.1.91'}
\langle\delta_{t_{m+1}}^-,f^{\star}_0\rangle=\mu([K^-,t_m[)\neq\mu([K^-,t_{m+1}[)\quad\Rightarrow\quad\langle\bm{\delta}_{\bb{t}}^-,f^{\star}\rangle\neq\langle\bm{\delta}_{\bb{t}}^-,f^{\star}_0\rangle.
\end{equation}
Since one can show that the value of $\langle\bm{\delta}_{\bb{t}}^-,\cdot\rangle$ over $\mathcal{V}(\K)$ is unique as in \citep[Theorem 3, Step 3 of the proof]{guillemet2025convergence}, \eqref{eq:4.1.91'} implies that $f_0^{\star}\notin\mathcal{V}(\K).$ This concludes the proof.
\end{proof}

\subsection{Proofs of Section 5.3}
\label{app:5.3}

\begin{proof}[\textbf{Proof of Theorem \ref{th:weakRT2}}]
The characterization of the solution set $\mathcal{V}(\K_{\epsilon})$ follows the proof of Theorem \ref{th:weakRT}. It remains to establish \eqref{eq:4.2.125}. 

Observe that $\mathcal{J}^{\star}(\K_{\epsilon})$ is decreasing in $\epsilon$ and is bounded from below. It is therefore sufficient to exhibit a sequence $(\mathcal{J}^{\star}(\K_{\epsilon_{\ell}}))_{\ell=1}^{\infty}$ such that
\begin{equation}
   \lim_{\ell\to\infty}\mathcal{J}^{\star}(\K_{\epsilon_{\ell}})=\mathcal{J}^{\star}(\K).
\end{equation}
Let $f_{\ell}=p_{\ell}+(u\ast \mu_{\ell})\in\M_{\Ls}(\K)$ be such that 
\begin{equation}
   \lim_{\ell\to\infty}\mathcal{J}(f_{\ell})=\mathcal{J}^{\star}(\K).
\end{equation}
We define 
\begin{equation}
    f_{\ell,\epsilon}=\mathrm{proj}_{\K_{\epsilon}}(f_{\ell})\in\M_{\Ls}(\K_{\epsilon}),
\end{equation}
where the projection operator is defined in the proof of Theorem \ref{th:weakRT}. Claims 1, 2, and 3 are as follows:  
\begin{itemize}
    \item [1.] the norm inequality $\norm{f_{\ell,\epsilon}}_{\M_{\Dd}}\leq\norm{f_{\ell}}_{\M_{\Dd}}$ holds;
    \item [2.] the measurement equality $\langle\bm{\delta}_{\bb{t}}^\omega,f_{\ell,\epsilon}\rangle=\langle\bm{\delta}_{\bb{t}}^\omega,f_{\ell}\rangle$ holds;
    \item [3.] the limit 
    $\underset{\epsilon\to\infty}{\mathrm{lim}\,}\langle\bm{\nu},f_{\ell,\epsilon}\rangle=\langle\bm{\nu},f_{\ell}\rangle$ holds.
\end{itemize}
Items 1 and 2 follow from the proof of Theorem \ref{th:weakRT}. Item 3 follows from the fact that $f_{\ell,\epsilon}$ is $\w$-convergent to $f_{\ell}$. Recall that  the lower-semicontinuous function $E(\bb{y},\cdot)$ is continuous on the interior of its domain, which here coincides with the entire domain. Consequently, there exists an $\epsilon_{\ell}$ sufficiently small such that 
\begin{equation}
    \mathcal{J}(f_{\ell,\epsilon_{\ell}})\leq\mathcal{J}(f_{\ell})+\frac{1}{\ell}.
\end{equation}
Therefore,
\begin{align}
        &\mathcal{J}^{\star}(\K)\leq\lim_{\ell\to\infty}\mathcal{J}(f_{\ell,\epsilon_{\ell}})\leq\lim_{\ell\to\infty}\mathcal{J}(f_{\ell})+\frac{1}{\ell}=\mathcal{J}^{\star}(\K)\nonumber\\
        \Rightarrow\,&\lim_{\ell\to\infty}\mathcal{J}(f_{\ell,\epsilon_{\ell}})=\mathcal{J}^{\star}(\K)\nonumber\\
        \Rightarrow\,&\lim_{\ell\to\infty}\mathcal{J}^{\star}(\K_{\epsilon_{\ell}})\leq\mathcal{J}^{\star}(\K)\nonumber\\
        \Rightarrow\,&\lim_{\ell\to\infty}\mathcal{J}^{\star}(\K_{\epsilon_{\ell}})=\mathcal{J}^{\star}(\K).
\end{align}
\end{proof}

\begin{proof}[\textbf{Proof of Theorem \ref{th:kuratowski}}]
We know from Theorem \ref{th:weakRT2} that $\mathcal{V}(\K_{\epsilon})$ is nonempty, convex, and $\w$-compact. Its extreme-point structure is also provided. The fact that $\bb{Ls}(\mathcal{V}(\K_{\epsilon}))$ is nonempty follows from the non-emptiness of each $\mathcal{V}(\K_{\epsilon})$ and their uniform boundedness. Its $\w$-compactness follows from the fact that $\bb{Ls}(\mathcal{V}(\K_{\epsilon}))$ is, by definition, closed, and can be shown to be bounded. 

Let $\tilde{f}\in\bb{Ls}(\mathcal{V}(\K_{\epsilon}))$. By definition, there exists a sequence $(f^{\star}_{\epsilon_{\ell}})_{\ell}^{\infty}$ that is $\w$-convergent to $\tilde f$, and the convergence $\underset{\epsilon\to0}{\mathrm{lim}}\mathcal{J}(f^{\star}_{\epsilon})=\mathcal{J}^{\star}(\K)$ follows from Theorem \ref{th:weakRT2}. It remains to derive the extreme-point structure of $\bb{Ls}(\mathcal{V}(\K_{\epsilon}))$. Observe that, for all $f_{\epsilon}\in\mathcal{V}(\K_{\epsilon})$, there exists a finite positive constant $C$ such that 
\begin{equation}
    \norm{f}_{\M_{\Dd}}\leq C.
\end{equation}
Define the set $E(\epsilon)$ as the set of extreme points of $\mathcal{V}(\K_{\epsilon})$, and 
\begin{equation}
    \overline E(\epsilon)=\Big\{f=c+\sum_{k=1}^{K}c_ku(\cdot-x_k):\quad\norm{f}_{\M_{\Dd}}\leq C,\,(\norm{c}_0+K)\leq M,\, x_k\in\K_{\epsilon}\Big\}.
\end{equation}
Observe that, $\forall\epsilon<\epsilon'$, one has that 
\begin{equation}
\label{eq:5.4.110}
    E(\epsilon)\subset\overline{E}(\epsilon)\quad\text{and}\quad\overline{E}(\epsilon')\subset\overline{E}(\epsilon)\subset\overline{E}(0).
\end{equation}
It follows from Choquet's theorem and \eqref{eq:5.4.110} that, for all $\ell$, there exists a probability measure $\rho_{\ell}\in\M(\overline{E}(0))$ such that $f^{\star}_{\epsilon_\ell}$ is represented as the barycentre
\begin{equation}
    f^{\star}_{\epsilon_\ell}=\int_{\overline{E}(0)}f\,\mathrm{d}\rho_{\ell}(f)
\end{equation}
in the sense that, for all $\phi$ in the predual  $\C_{\Dd}(\R)$, one has that 
\begin{equation}
    \langle f^{\star}_{\epsilon_\ell},\phi\rangle=\int_{\overline{E}(0)}\langle f,\phi\rangle\,\mathrm{d}\rho_{\ell}(f)=\int_{\mathcal{V}(\K_{\epsilon_{\ell}})\cap\overline{E}(0)}\langle f,\phi\rangle\,\mathrm{d}\rho_{\ell}(f).
\end{equation}
We observe that the $\w$-topology over $\overline{E}(0)$ is metrizable and therefore, that $\overline{E}(0)$ is a separable and complete (because it is compact) metric space. The Riesz-Markov-Kakutani shows that
\begin{equation}
\M(\overline{E}(0))=\C_0(\overline{E}(0))^{\star},\quad\text{where}\quad\rho_{\ell}\in\M(\overline{E}(0)).
\end{equation}
Then, the Banach-Alaoglu theorem reveals that there exists a subsequence (not relabeled) of $\{\rho_{\ell}\}_{\ell=1}^{\infty}$ that is $\w$-convergent to some limit $\tilde \rho$. We further notice that 
\begin{equation}
    \C_{\Dd}(\R)\subset\C_0(\overline{E}(0))
\end{equation}
such that $\forall\phi\in\C_\Dd(\R),\,\forall\rho\in\M(\overline{E}(0))$,
\begin{equation}
    \langle\rho,\phi\rangle=\int_{\overline{E}(0)}\langle f,\phi\rangle\,\mathrm{d}\rho(f).
\end{equation}
It follows directly that, $\forall\phi\in\C_\Dd(\R),$
\begin{equation}
    \langle \tilde f,\phi\rangle=\underset{\ell\to\infty}{\mathrm{lim}}\langle f^{\star}_{\epsilon_\ell},\phi\rangle=    \underset{\ell\to\infty}{\mathrm{lim}}\langle \rho_{\ell},\phi\rangle=\langle\tilde \rho,\phi\rangle=\int_{\bb{Ls}(\mathcal{V}(\K_{\epsilon}))\cap\overline E(0)}\langle f,\phi\rangle\,\mathrm{d}\rho(f),
\end{equation}
where the last equality follows from the fact that the probability measure $\rho_{\ell}$ is supported in $\mathcal{V}(\K_{\epsilon_{\ell}})\subset\overline E(0)$ and, because $\rho_{\ell}$ converges $\w$ to $\rho$, 
$\rho$ must therefore be supported in the limit superior of $\mathcal{V}(\K_{\epsilon_{\ell}})$, as given by $\bb{Ls}(\mathcal{V}(\K_{\epsilon_{\ell}}))$. In consequence, if $\tilde{f}$ is an extreme point of $\bb{Ls}(\mathcal{V}(\K_{\epsilon_{\ell}}))$, then the supporting measure $\rho$ must be a Dirac mass of the form $\delta_f$. This concludes the proof.
\end{proof}

\bibliographystyle{elsarticle-harv}
\bibliography{ref}

\end{document}